\documentclass[12pt,draftclsnofoot,onecolumn]{IEEEtran}

\usepackage{url,color}
\usepackage[dvips]{graphicx}
\usepackage{caption,subfig}
\usepackage{float}
\graphicspath{{../fig/}}
\DeclareGraphicsExtensions{.eps}
\usepackage{setspace}
\usepackage{amsmath}
\usepackage{amssymb}
\usepackage{cite}

\usepackage{array}
\usepackage{mdwmath}
\usepackage{mdwtab}

\begin{document}

\title{Random Distances Associated with Hexagons} 

\author{Yanyan Zhuang and Jianping Pan\\
University of Victoria, Victoria, BC, Canada}

\maketitle

\begin{abstract}
In this report, the explicit probability density functions of the random Euclidean distances
associated with regular hexagons are given, when the two endpoints of a link are randomly distributed in the same hexagon, and two adjacent
hexagons sharing a side, respectively. Simulation results show the accuracy of the obtained closed-form distance distribution
functions, which are important in a wide variety of
applied sciences and engineering fields. In particular, hexagons are often used
in wireless communication networks such as the cellular systems.
The correctness of these distance distribution functions is validated
by a recursion and a probabilistic sum. The first two statistical moments of the random distances,
and the polynomial fits of the density functions are also given in this report for practical uses. 
\end{abstract}

\begin{keywords}
Random distances; distance distribution functions; regular hexagons
\end{keywords}

\section{The Problem}

Define a ``unit hexagon'' as the regular hexagon with a side length of
$1$. Picking two points uniformly at random from the interior of a unit hexagon,
or between two adjacent unit hexagons sharing a side, the goal is to obtain the
explicit probability density function (PDF) of the random distances between these two endpoints. 

\section{Distance Distributions Associated with Regular Hexagons}\label{sec:result}

\subsection{Random Distances within a Unit Hexagon}

\begin{figure}
  \centering
  \includegraphics[width=0.4\columnwidth]{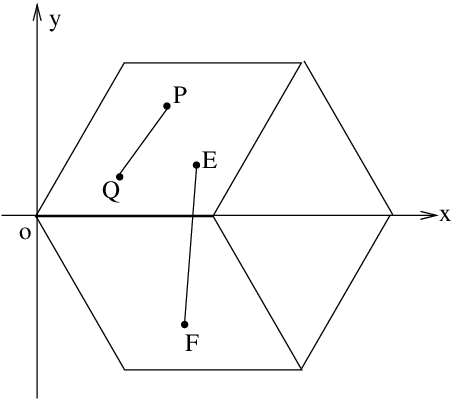}
  \caption{A Hexagon Decomposed into Three Rhombuses.}
  \label{fig:single_hexagon}
\end{figure}

A unit hexagon can be decomposed into three congruent rhombuses, as shown in
Fig.~\ref{fig:single_hexagon}. Define each one of these rhombuses as a ``unit rhombus'', i.e., with an acute angle of
$\frac{\pi}{3}$ and a side length of $1$. Here we first give the results of the random
distances associated with these rhombuses, and then use a probabilistic sum to combine them and obtain the
distribution of random distances for a unit hexagon.

\subsubsection{$|PQ|$}\label{sec:rhom1}

{\rm The probability density function of the random Euclidean distances
    between two uniformly distributed points that are both inside the same
    unit rhombus is\cite{rhombus}}
\begin{equation}\label{eq:fd_r_within}
  f_{D_1}(d)=2d\left\{
    \begin{array}{lr}

\left(\frac{4}{3}+\frac{2\pi}{9\sqrt{3}}\right)d^2-\frac{16}{3}d+\frac{2\pi}{
\sqrt{3}} & 0\leq d\leq \frac{\sqrt{3}}{2}\\

\frac{8}{\sqrt{3}}\left(1+\frac{d^2}{3}\right)\sin^{-1}\frac{\sqrt{3}}{2d}
+\left(\frac{4}{3}-\frac{10\pi}{9\sqrt{3}}\right)d^2-\frac{16}{3}d+\frac{10}{3}
\sqrt{4d^2-3}-\frac{2\pi}{\sqrt{3}} & \frac{\sqrt{3}}{2}\leq d\leq 1\\

\frac{4}{\sqrt{3}}\left(1-\frac{d^2}{3}\right)\sin^{-1}\frac{\sqrt{3}}{2d}
-\left(\frac{2}{3}-\frac{2\pi}{9\sqrt{3}}\right)d^2+\sqrt{4d^2-3}-\frac{2\pi}{
3\sqrt{3}}-1 & 1\leq d\leq \sqrt{3} \\

      0 & {\rm otherwise}
    \end{array}
  \right..
\end{equation}

\subsubsection{$|EF|$}\label{sec:rhom2}

{\rm The probability density function of the random Euclidean distances
    between two uniformly distributed points, one in each of the two adjacent 
    unit rhombuses sharing a side but with different orientation, as illustrated by ${\rm E}$ and ${\rm F}$ in Fig.~\ref{fig:single_hexagon}, is}
\begin{equation}\label{eq:fd_r_edge}
  f_{D_2}(d)=2d\left\{
    \begin{array}{lr}

\frac{4}{3}d-\left(\frac{1}{3}+\frac{2\pi}{9\sqrt{3}}\right)d^2 & 0\leq d\leq
\frac{\sqrt{3}}{2}\\

-\frac{4}{\sqrt{3}}\left(\frac{d^2}{3}+1\right)\sin^{-1}\frac{\sqrt{3}}{2d}
+\left(\frac{4\pi}{9\sqrt{3}}-\frac{1}{3}\right)d^2+\frac{4}{3}d-\frac{5}{3}
\sqrt{4d^2-3}+\frac{2\pi}{\sqrt{3}} & \frac{\sqrt{3}}{2}\leq d\leq 1\\

-\frac{2}{\sqrt{3}}\left(\frac{d^2}{3}+2\right)\sin^{-1}\frac{\sqrt{3}}{2d}
+\left(\frac{2\pi}{9\sqrt{3}}+\frac{1}{3}\right)d^2-\frac{3}{2}\sqrt{4d^2-3}
+\frac{2\pi}{\sqrt{3}}+\frac{1}{2} & 1\leq d\leq \sqrt{3} \\

\frac{2}{\sqrt{3}}\left(\frac{d^2}{3}+4\right)\sin^{-1}\frac{\sqrt{3}}{d}
-\left(\frac{2\pi}{9\sqrt{3}}+\frac{1}{3}\right)d^2+\frac{10}{3}\sqrt{d^2-3}
-\frac{8\pi}{3\sqrt{3}}-2 & \sqrt{3} \leq d\leq 2 \\

      0 & {\rm otherwise}
    \end{array}
  \right..
\end{equation}

\subsubsection{Final Result}\label{sec:hex1}

If we look at the two random endpoints of a given link inside a unit hexagon as shown in Fig.~\ref{fig:single_hexagon}, they will fall into
one of the two following cases: both endpoints are inside the same unit
rhombus, with probability $\frac{1}{3}$; each endpoint falls into one of the two adjacent
rhombuses sharing a side, with probability $\frac{2}{3}$. Therefore, 
given the results in Section~\ref{sec:rhom1} and \ref{sec:rhom2}, the probability density function of the random Euclidean distances
between these two endpoints is
$\frac{1}{3}f_{D_1}(d)+\frac{2}{3}f_{D_2}(d)$. We hence have the following:

\textbf{[Random distances within a unit hexagon]}
{\rm The probability density function of the random Euclidean distances
    between two uniformly distributed points that are both inside the same
    unit hexagon is}
\begin{equation}\label{eq:fd_h_within}
  f_{D_{\rm H_I}}(d)=\frac{2}{3}d\left\{
    \begin{array}{lr}

\left(\frac{2}{3}-\frac{2\pi}{9\sqrt{3}}\right)d^2-\frac{8}{3}d+\frac{2\pi}{
\sqrt{3}} & 0\leq d\leq 1\\

-\frac{4}{\sqrt{3}}\left(\frac{2d^2}{3}+1\right)\sin^{-1}\frac{\sqrt{3}}{2d}
+\frac{2\pi}{3\sqrt{3}}d^2-2\sqrt{4d^2-3}+\frac{10\pi}{3\sqrt{3}} &
1\leq d\leq \sqrt{3} \\

\frac{4}{\sqrt{3}}\left(\frac{d^2}{3}+4\right)\sin^{-1}\frac{\sqrt{3}}{d}
-\left(\frac{4\pi}{9\sqrt{3}}+\frac{2}{3}\right)d^2+\frac{20}{3}\sqrt{d^2-3}
-\frac{16\pi}{3\sqrt{3}}-4 & \sqrt{3} \leq d\leq 2 \\

      0 & {\rm otherwise}
    \end{array}
  \right..
\end{equation}

The corresponding cumulative distribution function (CDF) is 
\begin{equation}\label{eq:Fd_h_within}
  F_{D_{\rm H_I}}(d)=\left\{
    \begin{array}{lr}

\frac{1}{3}\left(\frac{1}{3}-\frac{\pi}{9\sqrt{3}}\right)d^4-\frac{16}{27}
d^3+\frac{2\pi}{3\sqrt{3}}d^2 & 0\leq d\leq 1\\

-\frac{4}{3\sqrt{3}}\left(\frac{d^4}{3}+d^2\right)\sin^{-1}\frac{\sqrt{3}}{2d}
+\frac{\pi}{9\sqrt{3}}d^4+\frac{10\pi}{9\sqrt{3}}d^2-\frac{26d^2+3}{54}\sqrt{
4d^2-3}+\frac{1}{18} & 1\leq d\leq \sqrt{3} \\

\frac{2}{3\sqrt{3}}\left(\frac{d^4}{3}+8d^2\right)\sin^{-1}\frac{\sqrt{3}}{d}
-\left(\frac{2\pi}{27\sqrt{3}}+\frac{1}{9}\right)d^4-\left(\frac{16\pi}{9\sqrt{3
}}+\frac{4}{3}\right)d^2 \\
~~~~~+\frac{14d^2+12}{9}\sqrt{d^2-3}+\frac{5}{9} & \sqrt{3} \leq d\leq 2 \\

      0 & {\rm otherwise}
    \end{array}
  \right..
\end{equation}

\subsection{Random Distances between Two Adjacent Unit Hexagons Sharing a Side}

\begin{figure}
  \centering
  \includegraphics[width=0.45\columnwidth]{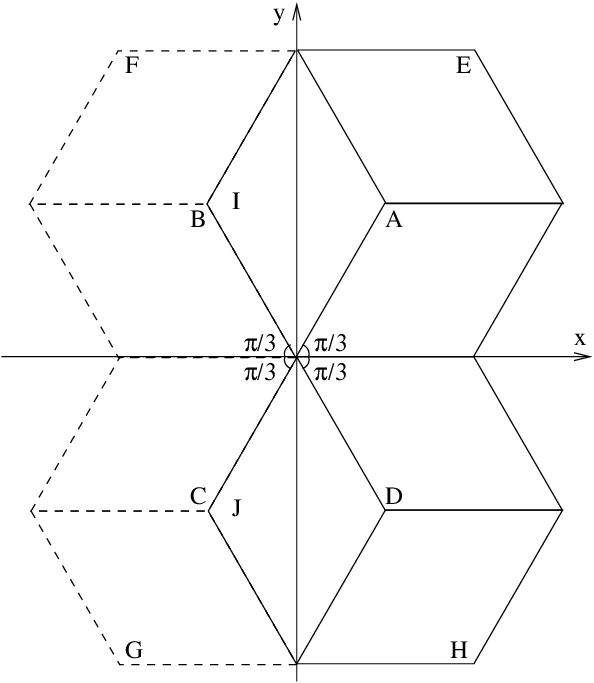}
  \caption{Random Points between Two Adjacent Hexagons: Different Cases with Rhombuses.}
  \label{fig:two_hexagon}
\end{figure}

Given two unit hexagons that are adjacent to each other, as shown in 
Fig.~\ref{fig:two_hexagon}, they can be decomposed into six congruent rhombuses 
in a similar way as that in Fig.~\ref{fig:single_hexagon}. Here we simply give the 
results for the random distances between these unit rhombuses with difference placement, as 
shown in Fig.~\ref{fig:two_hexagon}. At the end we combine them by a 
probabilistic sum, in order to obtain the distribution for the random distances between two 
adjacent hexagons sharing a side.

\subsubsection{$|DE|$}\label{sec:rhom3}

\begin{equation}\label{eq:fd_r_de}
  f_{D_3}(d)=2d\left\{
    \begin{array}{lr}

-\frac{2d^2}{3\sqrt{3}}\sin^{-1}\frac{\sqrt{3}}{2d}+\frac{\pi}{3\sqrt{3}}
d^2-\frac{\sqrt{4d^2-3}}{6} & \frac{\sqrt{3}}{2}\leq d\leq 1\\

\frac{d^2+4}{\sqrt{3}}\sin^{-1}\frac{\sqrt{3}}{2d}-\left(\frac{1}{3}+\frac{2\pi}
{9\sqrt{3}}\right)d^2-\frac{2}{3}d+\frac{19}{12}\sqrt{4d^2-3}-\frac{4\pi}{3\sqrt
{3}}-\frac{3}{4} & 1\leq d\leq \sqrt{3} \\

-\left(\frac{d^2}{3\sqrt{3}}+\frac{10}{\sqrt{3}}\right)\sin^{-1}\frac{\sqrt{3}}{
d}-\left(\frac{2d^2}{3\sqrt{3}}+2\sqrt{3}\right)\sin^{-1}\frac{\sqrt{3}}
{2d}+\left(\frac{4}{3}+\frac{2\pi}{9\sqrt{3}}\right)d^2 \\
~~~~~-\frac{13}{6}\sqrt{4d^2-3}-\frac{11}{3}\sqrt{d^2-3}-\frac{2}{3}d+\frac{
16\pi
}{3\sqrt{3}}+\frac{11}{2} & \sqrt{3} \leq d\leq 2 \\

\left(\frac{d^2}{3\sqrt{3}}+2\sqrt{3}\right)\sin^{-1}\frac{\sqrt{3}}{
d}-\left(\frac{2d^2}{3\sqrt{3}}+2\sqrt{3}\right)\sin^{-1}\frac{\sqrt{3}}
{2d}-\frac{13}{6}\sqrt{4d^2-3} \\
~~~~~+\frac{7}{3}\sqrt{d^2-3}+2d-\frac{1}{2} & 2\leq d \leq \frac{3\sqrt{3}}{2}
\\
 
\left(\frac{2d^2}{3\sqrt{3}}+4\sqrt{3}\right)\sin^{-1}\frac{3\sqrt{3}}{2d}
+\left(\frac{d^2}{3\sqrt{3}}+2\sqrt{3}\right)\sin^{-1}\frac{\sqrt{3}}{d} \\
~~~~-\left(\frac{2d^2}{3\sqrt{3}}+2\sqrt{3}\right)\sin^{-1}\frac{\sqrt{3}}
{2d}-\frac{\pi}{3\sqrt{3}}d^2-\frac{13}{6}\sqrt{4d^2-3}+\frac{7}{3}\sqrt{d^2-3}
& \frac{3\sqrt{3}}{2} \leq d \leq \sqrt{7} \\
~~~~+\frac{11}{6}\sqrt{4d^2-27}+2d-2\sqrt{3}\pi-\frac{1}{2}  \\

-\frac{d^2}{3\sqrt{3}}\sin^{-1}\frac{3\sqrt{3}}{2d}+\left(\frac{\pi}{9\sqrt{3}}
-\frac{1}{3}\right)d^2+2d-\frac{\sqrt{4d^2-27}}{4}-\frac{9}{4} & \sqrt{7} \leq d
\leq 3 \\

      0 & {\rm otherwise}
    \end{array}
  \right..
\end{equation}

\subsubsection{$|DI|$}\label{sec:rhom4}

\begin{equation}\label{eq:fd_r_di}
  f_{D_4}(d)=2d\left\{
    \begin{array}{lr}

\left(\frac{2\pi}{9\sqrt{3}}-\frac{1}{6}\right)d^2 & 0\leq d\leq
\frac{\sqrt{3}}{2}\\

\frac{2}{\sqrt{3}}\left(\frac{2d^2}{3}+1\right)\sin^{-1}\frac{\sqrt{3}}{
2d}-\left(\frac{4\pi}{9\sqrt{3}}+\frac{1}{6}\right)d^2+\sqrt{4d^2-3}-\frac{
\pi}{\sqrt{3}} & \frac{\sqrt{3}}{2}\leq d\leq 1 \\

\frac{1}{\sqrt{3}}\left(\frac{4d^2}{3}+1\right)\sin^{-1}\frac{\sqrt{3}}{
2d}+\left(\frac{1}{3}-\frac{4\pi}{9\sqrt{3}}\right)d^2+\frac{2}{3}\sqrt{
4d^2-3}-\frac{2}{3}d-\frac{2\pi}{3\sqrt{3}}+\frac{1}{2} & 1 \leq d\leq \sqrt{3}
\\

\frac{5}{\sqrt{3}}\sin^{-1}\frac{\sqrt{3}}{2d}-\frac{4}{\sqrt{3
}}\left(\frac{d^2}{3}+2\right)\sin^{-1}\frac{\sqrt{3}}{d}+\left(\frac{4\pi}
{9\sqrt{3}}-\frac{1}{3}\right)d^2+\frac{5}{3}\sqrt{4d^2-3} \\
~~~~-4\sqrt{d^2-3}-\frac{2}{3}d+\frac{8\pi}{3\sqrt{3}}-\frac{1}{2} &
\sqrt{3}\leq
d \leq 2 \\
 
\frac{5}{\sqrt{3}}\sin^{-1}\frac{\sqrt{3}}{2d}-\frac{2}{\sqrt{3
}}\left(\frac{d^2}{3}+2\right)\sin^{-1}\frac{\sqrt{3}}{d}+\left(\frac{2\pi}
{9\sqrt{3}}-\frac{1}{6}\right)d^2+\frac{5}{3}\sqrt{4d^2-3} \\
~~~~-2\sqrt{d^2-3}-2d+\frac{4\pi}{3\sqrt{3}}-\frac{1}{2} & 2 \leq d \leq
\sqrt{7}
\\

\left(\frac{2d^2}{3\sqrt{3}}+3\sqrt{3}\right)\sin^{-1}\frac{3\sqrt{3}}{
2d}+\left(\frac{1}{6}-\frac{2\pi}{9\sqrt{3}}\right)d^2+\frac{3}{2}\sqrt{4d^2-27}
-2d-\sqrt{3}\pi & \sqrt{7} \leq d \leq 3 \\

      0 & {\rm otherwise}
    \end{array}
  \right..
\end{equation}

\subsubsection{$|IJ|$}\label{sec:rhom5}

\begin{equation}\label{eq:ij}
  f_{D_5}(d)=2d\left\{
    \begin{array}{lr}

\left(\frac{1}{3}-\frac{\pi}{9\sqrt{3}}\right)d^2 & 0\leq d\leq 1\\

-\frac{4d^2}{3\sqrt{3}}\sin^{-1}\frac{\sqrt{3}}{2d}+\left(\frac{\pi}{3\sqrt{3}}
-1\right)d^2+\frac{8}{3}d-\frac{\sqrt{4d^2-3}}{3}-1 & 1\leq d\leq \sqrt{3}\\

\frac{4}{\sqrt{3}}\left(\frac{d^2}{3}-2\right)\sin^{-1}\frac{\sqrt{3
}}{2d}+\left(\frac{1}{3}-\frac{\pi}{9\sqrt{3}}\right)d^2+\frac{8}{3}d-\frac{7}{3
}\sqrt{4d^2-3} \\
~~~~+\frac{4\pi}{3\sqrt{3}}+1 & \sqrt{3}\leq d\leq 2\\

\frac{4}{\sqrt{3}}\left(\frac{d^2}{3}-2\right)\sin^{-1}\frac{\sqrt{3
}}{2d}+\frac{2d^2}{3\sqrt{3}}\sin^{-1}\frac{\sqrt{3}}{d}+\left(1-\frac{\pi}{
3\sqrt{3}}\right)d^2 \\
~~~~~-\frac{7}{3}\sqrt{4d^2-3}+\frac{2}{3}\sqrt{d^2-3}+\frac{4\pi}{3\sqrt{3}}+3
&
2\leq d\leq \sqrt{7} \\

\frac{2}{\sqrt{3}}\left(4-\frac{d^2}{3}\right)\sin^{-1}\frac{\sqrt{3}}{d}
+\left(\frac{\pi}{9\sqrt{3}}-\frac{1}{3}\right)d^2+2\sqrt{d^2-3} \\
~~~~-\frac{4\pi}{3\sqrt{3}}-2 & \sqrt{7} \leq d \leq 2\sqrt{3} \\

      0 & {\rm otherwise}
    \end{array}
  \right..
\end{equation}

\subsubsection{$|JE|$}\label{sec:rhom6}

\begin{equation}\label{eq:fd_r_je}
  f_{D_6}(d)=2d\left\{
    \begin{array}{lr}

-\left(\frac{d^2}{3\sqrt{3}}+\frac{1}{\sqrt{3}}\right)\sin^{-1}\frac{
\sqrt{3}}{2d}+\left(\frac{1}{6}+\frac{\pi}{9\sqrt{3}}\right)d^2-\frac{5}{12
}\sqrt{4d^2-3}+\frac{\pi}{3\sqrt{3}}+\frac{1}{4} & 1\leq d\leq \sqrt{3} \\

\left(\frac{2d^2}{3\sqrt{3}}+\frac{4}{\sqrt{3}}\right)\sin^{-1}\frac{\sqrt{3}}{d
}+\left(\frac{2}{\sqrt{3}}-\frac{d^2}{3\sqrt{3}}\right)\sin^{-1}\frac{\sqrt{3}}{
2d}-\left(\frac{1}{3}+\frac{2\pi}{9\sqrt{3}}\right)d^2 \\
~~~~+\frac{7}{12}\sqrt{4d^2-3}+2\sqrt{d^2-3}-\frac{13\pi}{6\sqrt{3}}-\frac{5}{4}
& \sqrt{3} \leq d\leq 2 \\

\frac{d^2}{3\sqrt{3}}\sin^{-1}\frac{\sqrt{3}}{d}+\left(\frac{2}{\sqrt{3}}-\frac{
d^2}{3\sqrt{3}}\right)\sin^{-1}\frac{\sqrt{3}}{2d}-\left(\frac{1}{6}
+\frac{\pi}{9\sqrt{3}}\right)d^2 \\
~~~~+\frac{7}{12}\sqrt{4d^2-3}+\frac{\sqrt{d^2-3}}{3}-\frac{5\pi}{6\sqrt{3}}
-\frac{1}{4} & 2 \leq d \leq \frac{3\sqrt{3}}{2} \\

\frac{d^2}{3\sqrt{3}}\sin^{-1}\frac{\sqrt{3}}{d}+\left(\frac{2}{\sqrt{3}
}-\frac{d^2}{3\sqrt{3}}\right)\sin^{-1}\frac{\sqrt{3}}{2d}-\left(\frac{2d^2}{
3\sqrt{3}}+3\sqrt{3}\right)\sin^{-1}\frac{3\sqrt{3}}{2d} \\
~~~~+\left(\frac{2\pi}{9\sqrt{3}}-\frac{1}{6}\right)d^2+\frac{7}{12}\sqrt{4d^2-3
}+\frac{\sqrt{d^2-3}}{3}-\frac{3}{2}\sqrt{4d^2-27} & \frac{3\sqrt{3}}{2} \leq d 
\leq \sqrt{7} \\
~~~~+\frac{11\pi}{3\sqrt{3}}-\frac{1}{4}  \\

\left(\frac{d^2}{3\sqrt{3}}-\frac{4}{\sqrt{3}}\right)\sin^{-1}\frac{\sqrt{3}}{d}
-\left(\frac{d^2}{3\sqrt{3}}+\frac{1}{2\sqrt{3}}\right)\sin^{-1}\frac{\sqrt{3}}{
2d} \\
~~~~-\left(\frac{2d^2}{3\sqrt{3}}+\frac{7\sqrt{3}}{2}\right)\sin^{-1}\frac{3\sqrt
{3}}{2d}+\left(\frac{1}{3}+\frac{2\pi}{9\sqrt{3}}\right)d^2-\frac{\sqrt{4d^2-3}}
{4} & \sqrt{7} \leq d \leq 3 \\
~~~~-\sqrt{d^2-3}-\frac{5}{3}\sqrt{4d^2-27}+\frac{11\pi}{2\sqrt{3}}+\frac{13}{4}
\\

\left(\frac{d^2}{3\sqrt{3}}-\frac{4}{\sqrt{3}}\right)\sin^{-1}\frac{\sqrt{3}}{d}
-\left(\frac{d^2}{3\sqrt{3}}+\frac{1}{2\sqrt{3}}\right)\sin^{-1}\frac{\sqrt{3}}{
2d}+\frac{5\sqrt{3}}{2}\sin^{-1}\frac{3\sqrt{3}}{2d} \\
~~~~-\frac{\sqrt{4d^2-3}}{4}-\sqrt{d^2-3}+\frac{5}{6}\sqrt{4d^2-27}-\frac{\pi}{
2\sqrt{3}}-\frac{5}{4}  & 3 \leq d \leq 2\sqrt{3} \\

\left(\frac{d^2}{3\sqrt{3}}+\frac{8}{\sqrt{3}}\right)\sin^{-1}\frac{2\sqrt{3}}{d
}-\left(\frac{d^2}{3\sqrt{3}}+\frac{1}{2\sqrt{3}}\right)\sin^{-1}\frac{\sqrt{3}}
{2d}+\frac{5\sqrt{3}}{2}\sin^{-1}\frac{3\sqrt{3}}{2d} \\
~~~~-\left(\frac{1}{6}+\frac{\pi}{9\sqrt{3}}\right)d^2-\frac{\sqrt{4d^2-3}}{4}
+\frac{5}{6}\sqrt{4d^2-27} +2\sqrt{d^2-12} & 2\sqrt{3} \leq d \leq \sqrt{13} \\
~~~~-\frac{31\pi}{6\sqrt{3}}-\frac{9}{4}  \\

      0 & {\rm otherwise}
    \end{array}
  \right..
\end{equation}

\subsubsection{$|HE|$}\label{sec:rhom7}

\begin{equation}\label{eq:fd_r_he}
  f_{D_7}(d)=2d\left\{
    \begin{array}{lr}

\frac{2}{\sqrt{3}}\left(\frac{2d^2}{3}+1\right)\sin^{-1}\frac{
\sqrt{3}}{2d}+\frac{4}{\sqrt{3}}\sin^{-1}\frac{\sqrt{3}}{d}-\left(\frac{1
}{3}+\frac{2\pi}{9\sqrt{3}}\right)d^2 \\
~~~~+\sqrt{4d^2-3}+\frac{4}{3}\sqrt{d^2-3}-\frac{7\pi}{3\sqrt{3}}-2 & \sqrt{3}
\leq d\leq 2 \\

\frac{2}{\sqrt{3}}\left(\frac{2d^2}{3}+1\right)\sin^{-1}\frac{
\sqrt{3}}{2d}-\frac{2}{\sqrt{3}}\left(\frac{d^2}{3}+2\right)\sin^{
-1}\frac{\sqrt{3}}{d} \\
~~~~+\sqrt{4d^2-3}-2\sqrt{d^2-3}+\frac{\pi}{3\sqrt{3}} & 2 \leq d \leq
\frac{3\sqrt{3}}{2} \\

\frac{2}{\sqrt{3}}\left(\frac{2d^2}{3}+1\right)\sin^{-1}\frac{
\sqrt{3}}{2d}-\frac{2}{\sqrt{3}}\left(\frac{d^2}{3}+2\right)\sin^{
-1}\frac{\sqrt{3}}{d}-2\sqrt{3}\sin^{-1}\frac{3\sqrt{3}}{2d} \\
~~~~+\sqrt{4d^2-3}-2\sqrt{d^2-3}-\frac{2}{3}\sqrt{4d^2-27}+\frac{10\pi}{3\sqrt{3
}
} & \frac{3\sqrt{3}}{2} \leq d \leq \sqrt{7} \\

\frac{1}{\sqrt{3}}\left(\frac{2d^2}{3}+1\right)\sin^{-1}\frac{
\sqrt{3}}{2d}-\frac{4}{\sqrt{3}}\left(\frac{d^2}{3}+2\right)\sin^{
-1}\frac{\sqrt{3}}{d}-3\sqrt{3}\sin^{-1}\frac{3\sqrt{3}}{2d} \\
~~~~+\left(\frac{1}{3}+\frac{2\pi}{9\sqrt{3}}\right)d^2+\frac{\sqrt{4d^2-3}}{2}
-4\sqrt{d^2-3}-\sqrt{ 4d^2-27} & \sqrt{7}\leq d \leq 3 \\
~~~~+\frac{17\pi}{3\sqrt{3}}+\frac{9}{2} \\

\frac{1}{\sqrt{3}}\left(\frac{2d^2}{3}+1\right)\sin^{-1}\frac{
\sqrt{3}}{2d}-\frac{4}{\sqrt{3}}\left(\frac{d^2}{3}+2\right)\sin^{
-1}\frac{\sqrt{3}}{d} \\
~~~~+\left(\frac{2d^2}{3\sqrt{3}}+3\sqrt{3}\right)\sin^{-1}\frac{3\sqrt{3}}{2d}
+\frac{\sqrt{4d^2-3}}{2}-4\sqrt{d^2-3} & 3 \leq d \leq 2\sqrt{3} \\
~~~~+\frac{3}{2}\sqrt{4d^2-27}-\frac{\pi}{3\sqrt{3}} \\

\frac{1}{\sqrt{3}}\left(\frac{2d^2}{3}+1\right)\sin^{-1}\frac{
\sqrt{3}}{2d}+\left(\frac{2d^2}{3\sqrt{3}}+3\sqrt{3}\right)\sin^{-1}\frac{
3\sqrt{3}}{2d}+\frac{8}{\sqrt{3}}\sin^{-1}\frac{2\sqrt{3}}{d} \\
~~~~-\left(\frac{1}{3}+\frac{2\pi}{9\sqrt{3}}\right)d^2+\frac{\sqrt{4d^2-3}}{2}
+\frac{3}{2}\sqrt{4d^2-27}+\frac{4}{3}\sqrt{d^2-12} & 2\sqrt{3} \leq d 
\leq \sqrt{13} \\
~~~~-\frac{17\pi}{3\sqrt{3}}-8 \\

      0 & {\rm otherwise}
    \end{array}
  \right..
\end{equation}

\subsubsection{Final Result}\label{sec:hex2}

Given the results in Section~\ref{sec:rhom1} to \ref{sec:rhom7}, when the two endpoints of a given link fall into one of the two adjacent 
hexagons sharing a side, the probability density function of the random distances between these two endpoints is
$\frac{1}{9}\left[f_{D_2}(d)+f_{D_5}(d)+f_{D_7}(d)\right]+\frac{2}{9}\left[f_{
D_3}(d)+f_{D_4}(d)+f_{D_6}(d)\right]$, using the similar reasoning as that in Section~\ref{sec:hex1}. With this probabilistic sum, we thus have the following:

\textbf{[Random distances between two unit hexagons]}
{\rm The probability density function of the random Euclidean distances
    between two uniformly distributed points, one in each of the two adjacent
    unit hexagons sharing a side, is}
\begin{equation}\label{eq:fd_h_btw}
  f_{D_{\rm H_A}}(d)=\frac{2}{9}d\left\{
    \begin{array}{lr}

\left(\frac{\pi}{9\sqrt{3}}-\frac{1}{3}\right)d^2+\frac{4}{3}d & 0\leq d\leq 1\\

\frac{2}{\sqrt{3}}(d^2+2)\sin^{-1}\frac{\sqrt{3}}{2d}-\left(\frac{1}{3}+\frac{
5\pi}{9\sqrt{3}}\right)d^2+\frac{11}{6}\sqrt{4d^2-3}-\frac{4\pi}{3\sqrt{3}}
-\frac{1}{2 } & 1\leq d\leq \sqrt{3} \\

\frac{2}{\sqrt{3}}\left(\frac{d^2}{3}-2\right)\sin^{-1}\frac{\sqrt{3}}{2d}-\frac
{4}{\sqrt{3}}\left(\frac{d^2}{3}+4\right)\sin^{-1}\frac{\sqrt{3}}{d}
+\left(1+\frac{\pi}{3\sqrt{3}}\right)d^2 \\
~~~~-\frac{7}{6}\sqrt{4d^2-3}-\frac{20}{3}\sqrt{d^2-3}+\frac{8\pi}{\sqrt{3}}
+\frac{9}{2} & \sqrt{3} \leq d\leq 2 \\

\frac{2}{\sqrt{3}}\left(\frac{d^2}{3}-2\right)\sin^{-1}\frac{\sqrt{3}}{2d}
+\left(\frac{1}{3}-\frac{\pi}{9\sqrt{3}}\right)d^2-\frac{7}{6}\sqrt{4d^2-3}
+\frac{8\pi}{3\sqrt{3}}+\frac{1}{2} & 2 \leq d \leq \sqrt{7} \\

-\frac{2}{\sqrt{3}}\left(\frac{d^2}{3}+6\right)\sin^{-1}\frac{3\sqrt{3}}{2d}
-\frac{4}{\sqrt{3}}\left(\frac{d^2}{3}+2\right)\sin^{-1}\frac{\sqrt{3}}{d}
+\left(\frac{1}{3}+\frac{5\pi}{9\sqrt{3}}\right)d^2 \\
~~~~-4\sqrt{d^2-3}-\frac{11}{6}\sqrt{4d^2-27}+\frac{28\pi}{3\sqrt{3}}+\frac{9}{2
}
& \sqrt{7} \leq d \leq 3 \\

\frac{2}{\sqrt{3}}\left(\frac{d^2}{3}+12\right)\sin^{-1}\frac{3\sqrt{3}}{2d}
-\frac{4}{\sqrt{3}}\left(\frac{d^2}{3}+2\right)\sin^{-1}\frac{\sqrt{3}}{d}
+\left(\frac{\pi}{9\sqrt{3}}-\frac{1}{3}\right)d^2 \\
~~~~-4\sqrt{d^2-3}+\frac{19}{6}\sqrt{4d^2-27}-\frac{8\pi}{3\sqrt{3}}-\frac{9}{2}
& 3 \leq d \leq 2\sqrt{3} \\

\frac{2}{\sqrt{3}}\left(\frac{d^2}{3}+12\right)\left(\sin^{-1}\frac{3\sqrt{3}}{
2d}+\sin^{-1}\frac{2\sqrt{3}}{d}\right)-\left(\frac{2}{3}+\frac{4\pi}{9\sqrt{3}}
\right)d^2\\
~~~~+\frac{19}{6}\sqrt{4d^2-27}+\frac{16}{3}\sqrt{d^2-12}-\frac{16\pi}{\sqrt{3}}
-\frac{25}{2} & 2\sqrt{3} \leq d\leq \sqrt{13} \\

      0 & {\rm otherwise}
    \end{array}
  \right..
\end{equation}

The corresponding CDF is 
\begin{equation}\label{eq:Fd_h_btw}
  F_{D_{\rm H_A}}(d)=\left\{
    \begin{array}{lr}

\frac{1}{18}\left(\frac{\pi}{9\sqrt{3}}-\frac{1}{3}\right)d^4+\frac{8}{81}d^3 &
0\leq d\leq 1\\

\frac{1}{9\sqrt{3}}\left(d^4+4d^2\right)\sin^{-1}\frac{\sqrt{3}}{2d}-\left(\frac
{5\pi}{162\sqrt{3}}+\frac{1}{54}\right)d^4-\left(\frac{4\pi}{27\sqrt{3}}+\frac{1
}{18}\right)d^2 \\
~~~~+\frac{94d^2+15}{648}\sqrt{4d^2-3}-\frac{1}{72} & 1\leq d\leq \sqrt{3} \\

\frac{1}{9\sqrt{3}}\left(\frac{d^4}{3}-4d^2\right)\sin^{-1}\frac{\sqrt{3}}{2d}
-\frac{2}{9\sqrt{3}}\left(\frac{d^4}{3}+8d^2\right)\sin^{-1}\frac{\sqrt{3}}{d}
\\
~~~~+\left(\frac{\pi}{54\sqrt{3}}+\frac{1}{18}\right)d^4+\left(\frac{8\pi}
{9\sqrt{3}}+\frac{1}{2}\right)d^2-\frac{2d^2+1}{24}\sqrt{4d^2-3} & \sqrt{3} 
\leq d\leq 2 \\
~~~~-\frac{4d^2+12}{27}\sqrt{d^2-3}-\frac{7}{72} \\

\frac{1}{9\sqrt{3}}\left(\frac{d^4}{3}-4d^2\right)\sin^{-1}\frac{\sqrt{3}}{2d}
-\left(\frac{\pi}{162\sqrt{3}}-\frac{1}{54}\right)d^4+\left(\frac{8\pi}{27\sqrt{
3 }}+\frac{1}{18}\right)d^2 \\
~~~~-\frac{2d^2+1}{24}\sqrt{4d^2-3}-\frac{53}{216} & 2 \leq d \leq \sqrt{7} \\

-\frac{2}{9\sqrt{3}}\left(\frac{d^4}{3}+4d^2\right)\sin^{-1}\frac{\sqrt{3}}{d}
-\frac{1}{3\sqrt{3}}\left(\frac{d^4}{9}+4d^2\right)\sin^{-1}\frac{3\sqrt{3}}{2d}
\\
~~~~+\left(\frac{5\pi}{162\sqrt{3}}+\frac{1}{54}\right)d^4+\left(\frac{28\pi}{
27\sqrt{3}}+\frac{1}{2}\right)d^2-\frac{26d^2+12}{81}\sqrt{d^2-3} & \sqrt{7}
\leq d \leq 3 \\
~~~~-\frac{94d^2+135}{648}\sqrt{4d^2-27}-\frac{101}{216} \\

-\frac{2}{9\sqrt{3}}\left(\frac{d^4}{3}+4d^2\right)\sin^{-1}\frac{\sqrt{3}}{d}
+\frac{1}{3\sqrt{3}}\left(\frac{d^4}{9}+8d^2\right)\sin^{-1}\frac{3\sqrt{3}}{2d}
\\
~~~~+\left(\frac{\pi}{162\sqrt{3}}-\frac{1}{54}\right)d^4-\left(\frac{8\pi}{
27\sqrt{3}}+\frac{1}{2}\right)d^2-\frac{26d^2+12}{81}\sqrt{d^2-3} & 3 \leq d
\leq 2\sqrt{3} \\
~~~~+\frac{158d^2+351}{648}\sqrt{4d^2-27}-\frac{263}{216} \\

\frac{1}{3\sqrt{3}}\left(\frac{d^4}{9}+8d^2\right)\left(\sin^{-1}\frac{3\sqrt{3}
}{2d}+\sin^{-1}\frac{2\sqrt{3}}{d}\right)-\left(\frac{2\pi}{81\sqrt{3}}+\frac{1}
{27}\right)d^4 \\
~~~-\left(\frac{16\pi}{9\sqrt{3}}+\frac{15}{18}\right)d^2+\frac{158d^2+351}{648}
\sqrt{4d^2-27} & 2\sqrt{3} \leq d\leq \sqrt{13} \\
~~~~+\frac{34d^2+96}{81}\sqrt{d^2-12}+\frac{25}{216} \\

      0 & {\rm otherwise}
    \end{array}
  \right..
\end{equation}

Note that although unit rhombuses and unit hexagons are assumed throughout 
(\ref{eq:fd_r_within})--(\ref{eq:Fd_h_btw}), the distance distribution
functions can be easily scaled by a nonzero scalar, for  rhombuses or hexagons of arbitrary 
side length. For example, let the side length of a regular hexagon be $s>0$, 
then
\begin{equation}
 F_{sD}(d)=P(sD\leq d)=P(D\leq \frac{d}{s})=F_D(\frac{d}{s}). \nonumber
\end{equation}
Therefore,
\begin{equation}\label{eq:scale}
 f_{sD}(d)=F'_D(\frac{d}{s})=\frac{1}{s}f_D(\frac{d}{s}),
\end{equation}
where $f_D(\cdot)$ can be (\ref{eq:fd_h_within}) or (\ref{eq:fd_h_btw}).

\section{Verification and Validation}
\subsection{Verification by Simulation}

\begin{figure}
  \centering
  \includegraphics[width=0.55\columnwidth]{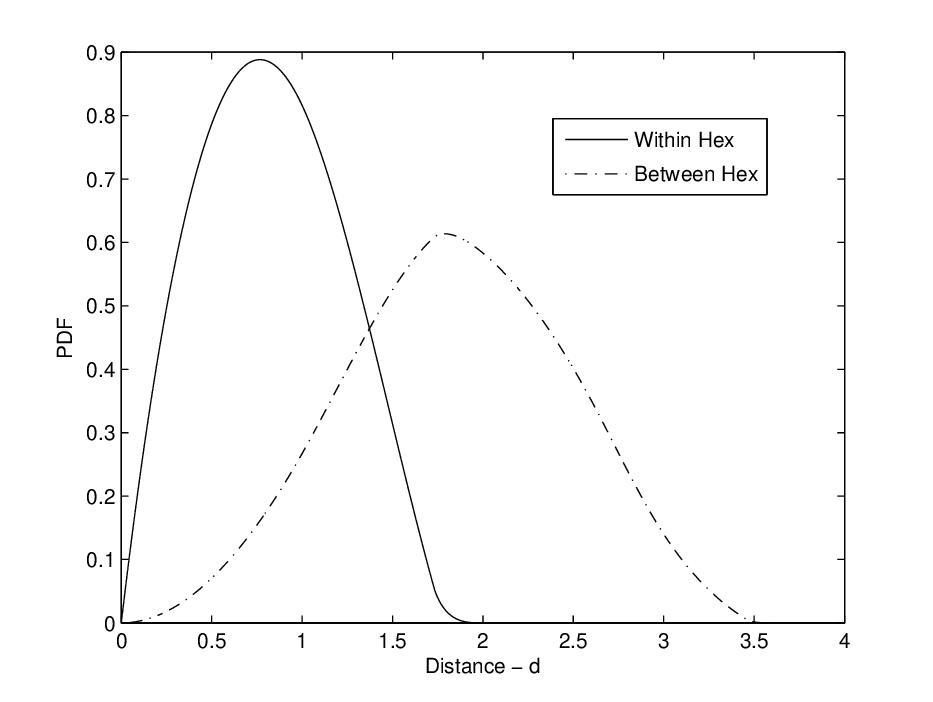}
  \caption{Distributions of Random Distances Associated with Hexagons.}
  \label{fig:hex_pdf}
\end{figure}

\begin{figure}
  \centering
  \includegraphics[width=0.55\columnwidth]{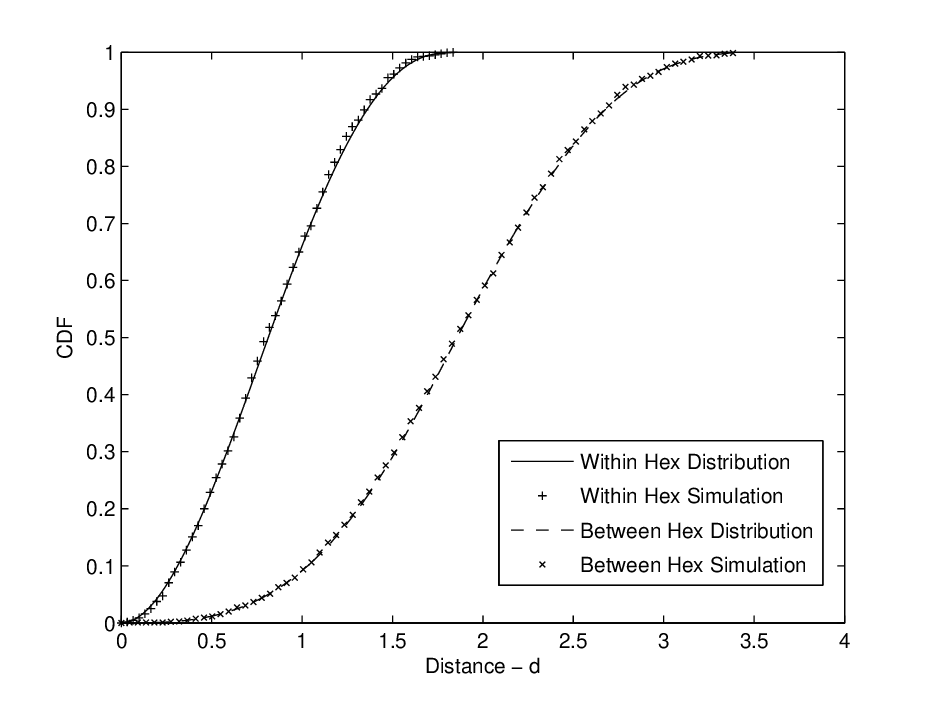}
  \caption{Distribution and Simulation Results For Random Distances Associated with
Hexagons.}
  \label{fig:hex_cdf}
\end{figure}

Figure~\ref{fig:hex_pdf} plots the probability density functions of the two
random distance cases given in (\ref{eq:fd_h_within}) and
(\ref{eq:fd_h_btw}), respectively. Figure~\ref{fig:hex_cdf} shows the comparison
between the cumulative distribution functions (CDFs) of the random distances,
and the simulation results by generating $2,000$ pairs of random points with the
corresponding geometric locations. Figure~\ref{fig:hex_cdf} demonstrates that
our distance distribution functions are very accurate when compared with the simulation results.

\subsection{Validation by Recursion}\label{sec:rec}

\begin{figure}
  \centering
  \includegraphics[width=0.35\columnwidth]{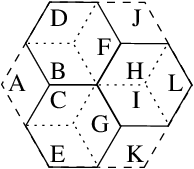}
  \caption{Partial Recursion through Hexagons and Rhombuses.}
  \label{fig:hexhex}
\end{figure}

As shown in Fig.~\ref{fig:hexhex}, a hexagon with a side length of $2$ can be
decomposed into three small hexagons with a side length of $1$, and three rhombuses ${\rm
A}$, ${\rm J}$ and ${\rm K}$, each with a side length of $1$ as well. With the scale transform in~(\ref{eq:scale}), the distance distribution in the large hexagon is 
$f_{2D}(d)=\frac{1}{2}f_{D_{\rm H_I}}(\frac{d}{2})$. On the other hand, if we look at the
two random endpoints of a given link inside the large hexagon, they  will fall into one of
the three following cases: i) both endpoints fall inside one of the
three small hexagons (${\rm BDF}$, ${\rm ECG}$ or ${\rm HIL}$), with
probability $\frac{3}{4}\times \frac{3}{4}$; ii) one of the endpoints falls
into one of the small hexagons, and the other endpoint falls into one of the three rhombuses ${\rm A}$,
${\rm J}$ or ${\rm K}$, with probability $2\times
\frac{3}{4}\times \frac{1}{4}$; iii) both endpoints fall into one of the rhombuses, with probability $\frac{1}{4}\times \frac{1}{4}$.

Each
of these three cases includes several more detailed sub-cases as follows:
\renewcommand{\labelenumi}{Case \roman{enumi})}
\begin{enumerate}
 \item 
%
Given the location of the first endpoint of a particular link, the second endpoint will fall in the same
hexagon as the first one with probability $\frac{1}{3}$, and in one of the two adjacent hexagons with probability
$\frac{2}{3}$. The unconditional
probability of these two sub-cases are $\frac{9}{16}\times
\frac{1}{3}=\frac{3}{16}$ and $\frac{9}{16}\times \frac{2}{3}=\frac{3}{8}$,
respectively.

 \item
%
Without loss of generality, suppose the first endpoint is located in ${\rm A}$. 
By symmetry, the second endpoint falls into any one of the four rhombuses ${\rm B}$, ${\rm
D}$, ${\rm F}$ and ${\rm H}$ with probability $\frac{2}{9}$, and into ${\rm L}$
with probability $\frac{1}{9}$. Thus the unconditional probabilities for $|AB|$,
$|AD|$, $|AF|$ or $|AH|$ are all $\frac{3}{8} \times \frac{2}{9}=\frac{1}{12}$, and
$\frac{3}{8} \times \frac{1}{9}=\frac{1}{24}$ for $|AL|$. 

 \item
%
If the first endpoint is in ${\rm A}$, then by symmetry, the second
endpoint is still located in ${\rm A}$ with probability $\frac{1}{3}$, and in
either one of ${\rm J}$ or ${\rm K}$ with probability $\frac{2}{3}$. The
unconditional probability of these two sub-cases are $\frac{1}{16}\times
\frac{1}{3}=\frac{1}{48}$ and $\frac{1}{16}\times \frac{2}{3}=\frac{1}{24}$,
respectively.

\end{enumerate}

In short, we have the probabilistic sum as
\begin{eqnarray}\label{eq:recur}
 f_{2D}(d)&=&\frac{3}{16}f_{D_{\rm H_I}}(d)+\frac{3}{8}f_{D_{\rm
H_A}}(d)+\frac{1}{12}\left[f_{D_2}(d)+f_{D_4}(d)+f_{D_3}(d)+f_{D_6}(d)\right]
\nonumber \\
&&+\frac{1}{24}\left[f_{D_7}(d)+f_{D_8}(d)\right]+\frac{1}{48}f_{D_1}(d),
\end{eqnarray}
where $f_{D_8}(d)$ is the density of $|AL|$ in Fig.~\ref{fig:hexhex}, which is the only 
distance distribution function that has not been given yet. We will do so in the immediate following.

\subsubsection{$|AL|$}

\begin{equation}
  f_{D_8}(d)=2d\left\{
    \begin{array}{lr}

-\frac{2}{\sqrt{3}}\left(\frac{d^2}{3}+4\right)\sin^{-1}\frac{\sqrt{3}}{d}
+\left(\frac{1}{3}+\frac{2\pi}{9\sqrt{3}}\right)d^2-\frac{10}{3}\sqrt{d^2-3}
+\frac{8\pi}{3\sqrt{3}}+2 & 2\leq d\leq \sqrt{7}\\

\frac{2}{\sqrt{3}}\left(\frac{d^2}{3}+8\right)\sin^{-1}\frac{\sqrt{3}}{d}
+\frac{4}{\sqrt{3}}\left(\frac{d^2}{3}+6\right)\sin^{-1}\frac{3\sqrt{3}}{2d}
-\left(1+\frac{2\pi}{3\sqrt{3}}\right)d^2 \\
~~~~+6\sqrt{d^2-3}+\frac{11}{3}\sqrt{4d^2-27}-\frac{40\pi}{3\sqrt{3}}-11 &
\sqrt{7}\leq d\leq 3\\

\frac{2}{\sqrt{3}}\left(\frac{d^2}{3}+8\right)\sin^{-1}\frac{\sqrt{3}}{d}
-\frac{4}{\sqrt{3}}\left(\frac{d^2}{3}+12\right)\sin^{-1}\frac{3\sqrt{3}}{2d}
+\left(\frac{1}{3}+\frac{2\pi}{9\sqrt{3}}\right)d^2 \\
~~~~+6\sqrt{d^2-3}-\frac{19}{3}\sqrt{4d^2-27}+\frac{32\pi}{3\sqrt{3}}+7 & 3\leq
d\leq 2\sqrt{3}\\

-\frac{4}{\sqrt{3}}\left(\frac{d^2}{3}+12\right)\sin^{-1}\frac{3\sqrt{3}}
{2d}-\frac{2}{\sqrt{3}}\left(\frac{d^2}{3}+8\right)\sin^{-1}\frac{2\sqrt{3}
}{d}+\left(1+\frac{2\pi}{3\sqrt{3}}\right)d^2 \\
~~~~-\frac{19}{3}\sqrt{4d^2-27}-4\sqrt{d^2-12}+\frac{64\pi}{3\sqrt{3}}+17 &
2\sqrt{3}\leq d\leq \sqrt{13} \\

\frac{2}{\sqrt{3}}\left(\frac{d^2}{3}+16\right)\sin^{-1}\frac{2\sqrt{3}
}{d}-\left(\frac{1}{3}+\frac{2\pi}{9\sqrt{3}}\right)d^2+\frac{20}{3}\sqrt{
d^2-12}-\frac{32\pi}{3\sqrt{3}}-8 & \sqrt{13} \leq d \leq 4 \\

      0 & {\rm otherwise}
    \end{array}
  \right..
\end{equation}

\subsubsection{Validation}

In order to confirm that the two definitions of $f_{2D}(d)$ at the beginning of Section~\ref{sec:rec} 
are equivalent, i.e., $\frac{3}{16}f_{D_{\rm H_I}}(d)+\frac{3}{8}f_{D_{\rm
H_A}}(d)+\frac{1}{12}\left[f_{D_2}(d)+f_{D_4}(d)+f_{D_3}(d)+f_{D_6}(d)\right]
+\frac{1}{24}\left[f_{D_7}(d)+f_{D_8}(d)\right]+\frac{1}{48}f_{D_1}(d)$ is equal to $\frac{1}{2}f_{D_{\rm H_I}}(\frac{d}{2})$,
we verify them mathematically as follows.

\renewcommand{\labelenumi}{\roman{enumi})}

\begin{enumerate}
 \item {$0\leq d\leq \frac{\sqrt{3}}{2}:$}
\begin{eqnarray}
\frac{3}{16}f_{D_{\rm
H_I}}(d)=\frac{d}{8}\left[\left(\frac{2}{3}-\frac{2\pi}{9\sqrt{3}}
\right)d^2-\frac{8}{3}d+\frac{2\pi}{\sqrt{3}}\right], \quad
\frac{3}{8}f_{D_{\rm
H_A}}(d)=\frac{d}{12}\left[\left(\frac{\pi}{9\sqrt{3}}-\frac{1}{3}
\right)d^2+\frac{4}{3}d\right], \nonumber
\end{eqnarray}
\begin{eqnarray}
\frac{1}{12}f_{D_2}(d)=\frac{d}{6}\left[\frac{4}{3}d-\left(\frac{1}{3}+\frac{
2\pi}{ 9\sqrt{3}}\right)d^2\right], \quad
\frac{1}{12}f_{D_4}(d)=\frac{d}{6}\left[\left(\frac{2\pi}{9\sqrt{3}}-\frac{1}{6}
\right)d^2\right], \nonumber
\end{eqnarray}
\begin{eqnarray}
\frac{1}{48}f_{D_1}(d)=\frac{d}{24}\left[\left(\frac{4}{3}+\frac{2\pi}{9\sqrt{3}
}\right)d^2-\frac{16}{3}d+\frac{2\pi}{\sqrt{3}}\right], \nonumber
\end{eqnarray}
while the probability density functions are $0$ for all other cases. Thus,
\begin{eqnarray}
f_{2D}(d)&=&\frac{3}{16}f_{D_{\rm H_I}}(d)+\frac{3}{8}f_{D_{\rm
H_A}}(d)+\frac{1}{12}\left[f_{D_2}(d)+f_{D_4}(d)\right]+\frac{1}{48}f_{D_1}(d)
\nonumber\\
&=&\frac{d}{6}\left[\left(\frac{1}{6}-\frac{\pi}{18\sqrt{3}}\right)d^2-\frac{4}{
3}d+\frac{2\pi}{\sqrt{3}}\right]=\frac{d}{6}\left[\left(\frac{2}{3}-\frac{2\pi}{
9\sqrt{3}}\right)\left(\frac{d}{2}\right)^2-\frac{8}{3}\left(\frac{d}{2}
\right)+\frac{2\pi}{\sqrt{3}}\right]\nonumber\\
&=&\frac{1}{2}f_{D_{\rm H_I}}(\frac{d}{ 2}). \nonumber
\end{eqnarray}

\item{$\frac{\sqrt{3}}{2} \leq d \leq 1:$}
\begin{eqnarray}
\frac{3}{16}f_{D_{\rm
H_I}}(d)=\frac{d}{8}\left[\left(\frac{2}{3}-\frac{2\pi}{9\sqrt{3}}
\right)d^2-\frac{8}{3}d+\frac{2\pi}{\sqrt{3}}\right], \quad
\frac{3}{8}f_{D_{\rm
H_A}}(d)=\frac{d}{12}\left[\left(\frac{\pi}{9\sqrt{3}}-\frac{1}{3}
\right)d^2+\frac{4}{3}d\right], \nonumber
\end{eqnarray}
\begin{eqnarray}
\frac{1}{12}f_{D_2}(d)=\frac{d}{6}\left[-\frac{4}{\sqrt{3}}\left(\frac{d^2}{3}
+1\right)\sin^{-1}\frac{\sqrt{3}}{2d}+\left(\frac{4\pi}{9\sqrt{3}}-\frac{1}{3}
\right)d^2+\frac{4}{3}d-\frac{5}{3}\sqrt{4d^2-3}+\frac{2\pi}{\sqrt{3}}\right],
\nonumber
\end{eqnarray}
\begin{eqnarray}
\frac{1}{12}f_{D_4}(d)=\frac{d}{6}\left[\frac{2}{\sqrt{3}}\left(\frac{2d^2}{3}
+1\right)\sin^{-1}\frac{\sqrt{3}}{2d}-\left(\frac{4\pi}{9\sqrt{3}}+\frac{1}{6}
\right)d^2+\sqrt{4d^2-3}-\frac{\pi}{\sqrt{3}}\right], \nonumber
\end{eqnarray}
\begin{eqnarray}
\frac{1}{12}f_{D_3}(d)=\frac{d}{6}\left[-\frac{2d^2}{3\sqrt{3}}\sin^{-1}\frac{
\sqrt{3}}{2d}+\frac{\pi}{3\sqrt{3}}d^2-\frac{\sqrt{4d^2-3}}{6}\right], \nonumber
\end{eqnarray}
\begin{eqnarray}
\frac{1}{48}f_{D_1}(d)=\frac{d}{24}\left[\frac{8}{\sqrt{3}}\left(1+\frac{d^2}{3}
\right)\sin^{-1}\frac{\sqrt{3}}{2d}+\left(\frac{4}{3}-\frac{10\pi}{9\sqrt{3}}
\right)d^2-\frac{16}{3}d+\frac{10}{3}\sqrt{4d^2-3}-\frac{2\pi}{\sqrt{3}}\right].
\nonumber
\end{eqnarray}
Thus,
\begin{eqnarray}
f_{2D}(d)&=&\frac{3}{16}f_{D_{\rm H_I}}(d)+\frac{3}{8}f_{D_{\rm
H_A}}(d)+\frac{1}{12}\left[f_{D_2}(d)+f_{D_4}(d)+f_{D_3}(d)\right]+\frac{1}{48}
f_{D_1}(d) \nonumber\\
&=&\frac{d}{6}\left[\left(\frac{2}{3}-\frac{2\pi}{9\sqrt{3}}\right)\left(\frac{
d}{2}\right)^2-\frac{8}{3}\left(\frac{d}{2}\right)+\frac{2\pi}{\sqrt{3}}\right]
=\frac{1}{2}f_{D_{\rm H_I}}(\frac{d}{2}). \nonumber
\end{eqnarray}

\item{$1\leq d \leq \sqrt{3}:$}
\begin{eqnarray}
\frac{3}{16}f_{D_{\rm
H_I}}(d)=\frac{d}{8}\left[-\frac{4}{\sqrt{3}}\left(\frac{2d^2}{3}+1\right)\sin^
{-1}\frac{\sqrt{3}}{2d}+\frac{2\pi}{3\sqrt{3}}d^2-2\sqrt{4d^2-3}+\frac{10\pi}{
3\sqrt{3}}\right], \nonumber
\end{eqnarray}
\begin{eqnarray}
\frac{3}{8}f_{D_{\rm
H_A}}(d)=\frac{d}{12}\left[\frac{2}{\sqrt{3}}(d^2+2)\sin^{-1}\frac{\sqrt{3}}{
2d}-\left(\frac{1}{3}+\frac{5\pi}{9\sqrt{3}}\right)d^2+\frac{11}{6}\sqrt{4d^2-3}
-\frac{4\pi}{3\sqrt{3}}-\frac{1}{2}\right], \nonumber
\end{eqnarray}
\begin{eqnarray}
\frac{1}{12}f_{D_2}(d)=\frac{d}{6}\left[-\frac{2}{\sqrt{3}}\left(\frac{d^2}{3}
+2\right)\sin^{-1}\frac{\sqrt{3}}{2d}+\left(\frac{2\pi}{9\sqrt{3}}+\frac{1}{3}
\right)d^2-\frac{3}{2}\sqrt{4d^2-3}+\frac{2\pi}{\sqrt{3}}+\frac{1}{2}\right],
\nonumber
\end{eqnarray}
\begin{eqnarray}
\frac{1}{12}f_{D_4}(d)=\frac{d}{6}\left[\frac{1}{\sqrt{3}}\left(\frac{4d^2}{3}
+1\right)\sin^{-1}\frac{\sqrt{3}}{2d}+\left(\frac{1}{3}-\frac{4\pi}{9\sqrt{3}}
\right)d^2+\frac{2}{3}\sqrt{4d^2-3}-\frac{2}{3}d-\frac{2\pi}{3\sqrt{3}}+\frac{1}
{2}\right], \nonumber
\end{eqnarray}
\begin{eqnarray}
\frac{1}{12}f_{D_3}(d)=\frac{d}{6}\left[\frac{d^2+4}{\sqrt{3}}\sin^{-1}\frac{
\sqrt{3}
}{2d}-\left(\frac{1}{3}+\frac{2\pi}{9\sqrt{3}}\right)d^2-\frac{2}{3}d+\frac{19}{
12}\sqrt{4d^2-3}-\frac{4\pi}{3\sqrt{3}}-\frac{3}{4}\right], \nonumber
\end{eqnarray}
\begin{eqnarray}
\frac{1}{12}f_{D_6}(d)=\frac{d}{6}\left[-\left(\frac{d^2}{3\sqrt{3}}+\frac{1}{
\sqrt{3}}\right)\sin^{-1}\frac{\sqrt{3}}{2d}+\left(\frac{1}{6}+\frac{\pi}{9\sqrt
{3}}\right)d^2-\frac{5}{12}\sqrt{4d^2-3}+\frac{\pi}{3\sqrt{3}}+\frac{1}{4}\right
], \nonumber
\end{eqnarray}
\begin{eqnarray}
\frac{1}{48}f_{D_1}(d)=\frac{d}{24}\left[\frac{4}{\sqrt{3}}\left(1-\frac{d^2}{3}
\right)\sin^{-1}\frac{\sqrt{3}}{2d}-\left(\frac{2}{3}-\frac{2\pi}{9\sqrt{3}}
\right)d^2+\sqrt{4d^2-3}-\frac{2\pi}{3\sqrt{3}}-1\right]. \nonumber
\end{eqnarray}
Thus,
\begin{eqnarray}
f_{2D}(d)&=&\frac{3}{16}f_{D_{\rm H_I}}(d)+\frac{3}{8}f_{D_{\rm
H_A}}(d)+\frac{1}{12}\left[f_{D_2}(d)+f_{D_4}(d)+f_{D_3}(d)+f_{D_6}(d)\right]
\nonumber\\
&&+\frac{1}{48}f_{D_1}(d)=\frac{d}{6}\left[\left(\frac{2}{3}-\frac{2\pi}{9\sqrt{
3}}\right)\left(\frac{d}{2}\right)^2-\frac{8}{3}\left(\frac{d}{2}\right)+\frac{
2\pi}{\sqrt{3}}\right]=\frac{1}{2}f_{D_{\rm H_I}}(\frac{d}{2}). \nonumber
\end{eqnarray}

\item{$\sqrt{3} \leq d \leq 2:$}
\begin{eqnarray}
\frac{3}{16}f_{D_{\rm
H_I}}(d)=\frac{d}{8}\left[\frac{4}{\sqrt{3}}\left(\frac{d^2}{3}+4\right)\sin^{
-1}\frac{\sqrt{3}}{d}-\left(\frac{4\pi}{9\sqrt{3}}+\frac{2}{3}\right)d^2+\frac{
20}{3}\sqrt{d^2-3}-\frac{16\pi}{3\sqrt{3}}-4\right], \nonumber
\end{eqnarray}
\begin{eqnarray}
\frac{3}{8}f_{D_{\rm
H_A}}(d)&=&\frac{d}{12}\left[\frac{2}{\sqrt{3}}\left(\frac{d^2}{3}
-2\right)\sin^{-1}\frac{\sqrt{3}}{2d}-\frac{4}{\sqrt{3}}\left(\frac{d^2}{3}
+4\right)\sin^{-1}\frac{\sqrt{3}}{d} \right. \nonumber\\
&&\left.+\left(1+\frac{\pi}{3\sqrt{3}}\right)d^2-\frac{7}{6}\sqrt{4d^2-3}-\frac{
20}{3}\sqrt{d^2-3}+\frac{8\pi}{\sqrt{3}}+\frac{9}{2}\right], \nonumber
\end{eqnarray}
\begin{eqnarray}
\frac{1}{12}f_{D_2}(d)=\frac{d}{6}\left[\frac{2}{\sqrt{3}}\left(\frac{d^2}{3}
+4\right)\sin^{-1}\frac{\sqrt{3}}{d}-\left(\frac{2\pi}{9\sqrt{3}}+\frac{1}{3}
\right)d^2+\frac{10}{3}\sqrt{d^2-3}-\frac{8\pi}{3\sqrt{3}}-2\right], \nonumber
\end{eqnarray}
\begin{eqnarray}
\frac{1}{12}f_{D_4}(d)&=&\frac{d}{6}\left[\frac{5}{\sqrt{3}}\sin^{-1}\frac{\sqrt
{3}}{2d}-\frac{4}{\sqrt{3}}\left(\frac{d^2}{3}+2\right)\sin^{-1}\frac{\sqrt{3}}{
d}+\left(\frac{4\pi}{9\sqrt{3}}-\frac{1}{3}\right)d^2+\frac{5}{3}\sqrt{4d^2-3}
\right. \nonumber\\
&&\left.-4\sqrt{d^2-3}-\frac{2}{3}d+\frac{8\pi}{3\sqrt{3}}-\frac{1}{2}\right],
\nonumber
\end{eqnarray}
\begin{eqnarray}
\frac{1}{12}f_{D_3}(d)&=&\frac{d}{6}\left[-\left(\frac{d^2}{3\sqrt{3}}+\frac{10}
{\sqrt{3}}\right)\sin^{-1}\frac{\sqrt{3}}{d}-\left(\frac{2d^2}{3\sqrt{3}}+2\sqrt
{3}\right)\sin^{-1}\frac{\sqrt{3}}{2d}+\left(\frac{4}{3}+\frac{2\pi}{9\sqrt{3}}
\right)d^2\right. \nonumber\\
&&\left.-\frac{13}{6}\sqrt{4d^2-3}-\frac{11}{3}\sqrt{d^2-3}-\frac{2}{3}d+\frac{
16\pi}{3\sqrt{3}}+\frac{11}{2}\right], \nonumber
\end{eqnarray}
\begin{eqnarray}
\frac{1}{12}f_{D_6}(d)&=&\frac{d}{6}\left[\left(\frac{2d^2}{3\sqrt{3}}+\frac{4}{
\sqrt{3}}\right)\sin^{-1}\frac{\sqrt{3}}{d}+\left(\frac{2}{\sqrt{3}}-\frac{d^2}{
3\sqrt{3}}\right)\sin^{-1}\frac{\sqrt{3}}{2d}-\left(\frac{1}{3}+\frac{2\pi}{
9\sqrt{3}}\right)d^2\right. \nonumber\\
&&\left.+\frac{7}{12}\sqrt{4d^2-3}+2\sqrt{d^2-3}-\frac{13\pi}{6\sqrt{3}}-\frac{5
}{4}\right], \nonumber
\end{eqnarray}
\begin{eqnarray}
\frac{1}{24}f_{D_7}(d)&=&\frac{d}{12}\left[\frac{2}{\sqrt{3}}\left(\frac{2d^2}{3
}+1\right)\sin^{-1}\frac{\sqrt{3}}{2d}+\frac{4}{\sqrt{3}}\sin^{-1}\frac{\sqrt{3}
}{d}-\left(\frac{1}{3}+\frac{2\pi}{9\sqrt{3}}\right)d^2+\sqrt{4d^2-3}\right.
\nonumber\\
&&\left.+\frac{4}{3}\sqrt{d^2-3}-\frac{7\pi}{3\sqrt{3}}-2\right]. \nonumber
\end{eqnarray}
Thus,
\begin{eqnarray}
f_{2D}(d)&=&\frac{3}{16}f_{D_{\rm H_I}}(d)+\frac{3}{8}f_{D_{\rm
H_A}}(d)+\frac{1}{12}\left[f_{D_2}(d)+f_{D_4}(d)+f_{D_3}(d)+f_{D_6}(d)\right]
\nonumber\\
&&+\frac{1}{24}f_{D_7}(d)=\frac{d}{6}\left[\left(\frac{2}{3}-\frac{2\pi}{9\sqrt{
3}}\right)\left(\frac{d}{2}\right)^2-\frac{8}{3}\left(\frac{d}{2}\right)+\frac{
2\pi}{\sqrt{3}}\right]=\frac{1}{2}f_{D_{\rm H_I}}(\frac{d}{2}). \nonumber
\end{eqnarray}

\item{$2\leq d \leq \frac{3\sqrt{3}}{2}:$}
\begin{eqnarray}
\frac{3}{8}f_{D_{\rm
H_A}}(d)=\frac{d}{12}\left[\frac{2}{\sqrt{3}}\left(\frac{d^2}{3}-2\right)\sin^
{-1}\frac{\sqrt{3}}{2d}+\left(\frac{1}{3}-\frac{\pi}{9\sqrt{3}}\right)d^2-\frac{
7}{6}\sqrt{4d^2-3}+\frac{8\pi}{3\sqrt{3}}+\frac{1}{2}\right], \nonumber
\end{eqnarray}
\begin{eqnarray}
\frac{1}{12}f_{D_4}(d)&=&\frac{d}{6}\left[\frac{5}{\sqrt{3}}\sin^{-1}\frac{\sqrt
{3}}{2d}-\frac{2}{\sqrt{3}}\left(\frac{d^2}{3}+2\right)\sin^{-1}\frac{\sqrt{3}}{
d}+\left(\frac{2\pi}{9\sqrt{3}}-\frac{1}{6}\right)d^2+\frac{5}{3}\sqrt{4d^2-3}
\right. \nonumber\\
&&\left.-2\sqrt{d^2-3}-2d+\frac{4\pi}{3\sqrt{3}}-\frac{1}{2}\right], \nonumber
\end{eqnarray}
\begin{eqnarray}
\frac{1}{12}f_{D_3}(d)&=&\frac{d}{6}\left[\left(\frac{d^2}{3\sqrt{3}}+2\sqrt{3}
\right)\sin^{-1}\frac{\sqrt{3}}{d}-\left(\frac{2d^2}{3\sqrt{3}}+2\sqrt{3}
\right)\sin^{-1}\frac{\sqrt{3}}{2d}-\frac{13}{6}\sqrt{4d^2-3}\right. \nonumber\\
&&\left.+\frac{7}{3}\sqrt{d^2-3}+2d-\frac{1}{2}\right], \nonumber
\end{eqnarray}
\begin{eqnarray}
\frac{1}{12}f_{D_6}(d)&=&\frac{d}{6}\left[\frac{d^2}{3\sqrt{3}}\sin^{-1}\frac{
\sqrt{3}}{d}+\left(\frac{2}{\sqrt{3}}-\frac{d^2}{3\sqrt{3}}\right)\sin^{-1}\frac
{\sqrt{3}}{2d}-\left(\frac{1}{6}+\frac{\pi}{9\sqrt{3}}\right)d^2+\frac{7}{12}
\sqrt{4d^2-3}\right. \nonumber\\
&&\left.+\frac{\sqrt{d^2-3}}{3}-\frac{5\pi}{6\sqrt{3}}-\frac{1}{4}\right],
\nonumber
\end{eqnarray}
\begin{eqnarray}
\frac{1}{24}f_{D_7}(d)&=&\frac{d}{12}\left[\frac{2}{\sqrt{3}}\left(\frac{2d^2}{3
}
+1\right)\sin^{-1}\frac{\sqrt{3}}{2d}-\frac{2}{\sqrt{3}}\left(\frac{d^2}{3}
+2\right)\sin^{-1}\frac{\sqrt{3}}{d}+\sqrt{4d^2-3}\right. \nonumber \\
&&\left.-2\sqrt{d^2-3}+\frac{\pi}{3\sqrt{3}}\right], \nonumber
\end{eqnarray}
\begin{eqnarray}
\frac{1}{24}f_{D_8}(d)&=&\frac{d}{12}\left[-\frac{2}{\sqrt{3}}\left(\frac{d^2}{3
}+4\right)\sin^{-1}\frac{\sqrt{3}}{d}+\left(\frac{1}{3}+\frac{2\pi}{9\sqrt{3}}
\right)d^2-\frac{10}{3}\sqrt{d^2-3}+\frac{8\pi}{3\sqrt{3}}+2\right]. \nonumber
\end{eqnarray}
Thus,
\begin{eqnarray}
f_{2D}(d)&=&\frac{3}{8}f_{D_{\rm
H_A}}(d)+\frac{1}{12}\left[f_{D_4}(d)+f_{D_3}(d)+f_{D_6}(d)\right]+\frac{1}{24}
\left[f_{D_7}(d)+f_{D_8}(d)\right] \nonumber\\
&=&\frac{d}{6}\left[-\frac{2}{\sqrt{3}}\left(\frac{d^2}{3}+2\right)\sin^{-1}
\frac{\sqrt{3}}{d}+\frac{\pi}{6\sqrt{3}}d^2-2\sqrt{d^2-3}+\frac{10\pi}{3\sqrt{3}
}\right] \nonumber\\
&=&\frac{d}{6}\left[-\frac{4}{\sqrt{3}}\left(\frac{2(d/2)^2}{3}+1\right)\sin^{-1
}\frac{\sqrt{3}}{2(d/2)}+\frac{2\pi}{3\sqrt{3}}\left(\frac{d}{2}\right)^2-2\sqrt
{4\left(\frac{d}{2}\right)^2-3}+\frac{10\pi}{3\sqrt{3}}\right]\nonumber \\
&=&\frac{1}{2}f_{D_{\rm H_I}}(\frac{d}{2}). \nonumber
\end{eqnarray}

\item{$\frac{3\sqrt{3}}{2} \leq d\leq \sqrt{7}:$}
\begin{eqnarray}
\frac{3}{8}f_{D_{\rm
H_A}}(d)=\frac{d}{12}\left[\frac{2}{\sqrt{3}}\left(\frac{d^2}{3}-2\right)\sin^
{-1}\frac{\sqrt{3}}{2d}+\left(\frac{1}{3}-\frac{\pi}{9\sqrt{3}}\right)d^2-\frac{
7}{6}\sqrt{4d^2-3}+\frac{8\pi}{3\sqrt{3}}+\frac{1}{2}\right], \nonumber
\end{eqnarray}
\begin{eqnarray}
\frac{1}{12}f_{D_4}(d)&=&\frac{d}{6}\left[\frac{5}{\sqrt{3}}\sin^{-1}\frac{\sqrt
{3}}{2d}-\frac{2}{\sqrt{3}}\left(\frac{d^2}{3}+2\right)\sin^{-1}\frac{\sqrt{3}}{
d}+\left(\frac{2\pi}{9\sqrt{3}}-\frac{1}{6}\right)d^2+\frac{5}{3}\sqrt{4d^2-3}
\right. \nonumber\\
&&\left.-2\sqrt{d^2-3}-2d+\frac{4\pi}{3\sqrt{3}}-\frac{1}{2}\right], \nonumber
\end{eqnarray}
\begin{eqnarray}
\frac{1}{12}f_{D_3}(d)&=&\frac{d}{6}\left[\left(\frac{2d^2}{3\sqrt{3}}+4\sqrt{3}
\right)\sin^{-1}\frac{3\sqrt{3}}{2d}+\left(\frac{d^2}{3\sqrt{3}}+2\sqrt{3}
\right)\sin^{-1}\frac{\sqrt{3}}{d} \right. \nonumber\\
&&\left.-\left(\frac{2d^2}{3\sqrt{3}}+2\sqrt{3}
\right)\sin^{-1}\frac{\sqrt{3}}{2d}-\frac{\pi}{3\sqrt{3}}d^2-\frac{13}{6}\sqrt{
4d^2-3}+\frac{7}{3}\sqrt{d^2-3}\right. \nonumber\\
&&\left.+\frac{11}{6}\sqrt{4d^2-27}+2d-2\sqrt{3}\pi-\frac{1}{2}\right],
\nonumber
\end{eqnarray}
\begin{eqnarray}
\frac{1}{12}f_{D_6}(d)&=&\frac{d}{6}\left[\frac{d^2}{3\sqrt{3}}\sin^{-1}\frac{
\sqrt{3}}{d}+\left(\frac{2}{\sqrt{3}}-\frac{d^2}{3\sqrt{3}}\right)\sin^{-1}\frac
{\sqrt{3}}{2d}-\left(\frac{2d^2}{3\sqrt{3}}+3\sqrt{3}\right)\sin^{-1}\frac{
3\sqrt{3}}{2d}\right. \nonumber\\
&&\left.+\left(\frac{2\pi}{9\sqrt{3}}-\frac{1}{6}\right)d^2+\frac{7}{12}\sqrt{
4d^2-3}+\frac{\sqrt{d^2-3}}{3}-\frac{3}{2}\sqrt{4d^2-27}+\frac{11\pi}{3\sqrt{3}}
-\frac{1}{4}\right], \nonumber
\end{eqnarray}
\begin{eqnarray}
\frac{1}{24}f_{D_7}(d)&=&\frac{d}{12}\left[\frac{2}{\sqrt{3}}\left(\frac{2d^2}{3
}+1\right)\sin^{-1}\frac{\sqrt{3}}{2d}-\frac{2}{\sqrt{3}}\left(\frac{d^2}{3}
+2\right)\sin^{-1}\frac{\sqrt{3}}{d}-2\sqrt{3}\sin^{-1}\frac{3\sqrt{3}}{2d}
\right. \nonumber\\
&&\left.+\sqrt{4d^2-3}-2\sqrt{d^2-3}-\frac{2}{3}\sqrt{4d^2-27}+\frac{10\pi}{
3\sqrt{3}}\right], \nonumber
\end{eqnarray}
\begin{eqnarray}
\frac{1}{24}f_{D_8}(d)&=&\frac{d}{12}\left[-\frac{2}{\sqrt{3}}\left(\frac{d^2}{3
}+4\right)\sin^{-1}\frac{\sqrt{3}}{d}+\left(\frac{1}{3}+\frac{2\pi}{9\sqrt{3}}
\right)d^2-\frac{10}{3}\sqrt{d^2-3}+\frac{8\pi}{3\sqrt{3}}+2\right]. \nonumber
\end{eqnarray}
Thus,
\begin{eqnarray}
f_{2D}(d)&=&\frac{3}{8}f_{D_{\rm
H_A}}(d)+\frac{1}{12}\left[f_{D_4}(d)+f_{D_3}(d)+f_{D_6}(d)\right]+\frac{1}{24}
\left[f_{D_7}(d)+f_{D_8}(d)\right] \nonumber\\
&=&\frac{d}{6}\left[-\frac{4}{\sqrt{3}}\left(\frac{2(d/2)^2}{3}+1\right)\sin^{-1
}\frac{\sqrt{3}}{2(d/2)}+\frac{2\pi}{3\sqrt{3}}\left(\frac{d}{2}\right)^2-2\sqrt
{4\left(\frac{d}{2}\right)^2-3}+\frac{10\pi}{3\sqrt{3}}\right]\nonumber \\
&=&\frac{1}{2}f_{D_{\rm H_I}}(\frac{d}{2}). \nonumber
\end{eqnarray}

\item{$\sqrt{7} \leq d\leq 3:$}
\begin{eqnarray}
\frac{3}{8}f_{D_{\rm
H_A}}(d)&=&\frac{d}{12}\left[-\frac{2}{\sqrt{3}}\left(\frac{d^2}{3}
+6\right)\sin^{-1}\frac{3\sqrt{3}}{2d}-\frac{4}{\sqrt{3}}\left(\frac{d^2}{3}
+2\right)\sin^{-1}\frac{\sqrt{3}}{d} \right. \nonumber\\
&&\left.+\left(\frac{1}{3}+\frac{5\pi}{9\sqrt{3}}\right)d^2-4\sqrt{d^2-3}-\frac{
11}{6}\sqrt{4d^2-27}+\frac{28\pi}{3\sqrt{3}}+\frac{9}{2}\right], \nonumber
\end{eqnarray}
\begin{eqnarray}
\frac{1}{12}f_{D_4}(d)&=&\frac{d}{6}\left[\left(\frac{2d^2}{3\sqrt{3}}+3\sqrt{3}
\right)\sin^{-1}\frac{3\sqrt{3}}{2d}+\left(\frac{1}{6}-\frac{2\pi}{9\sqrt{3}}
\right)d^2+\frac{3}{2}\sqrt{4d^2-27}-2d-\sqrt{3}\pi\right], \nonumber
\end{eqnarray}
\begin{eqnarray}
\frac{1}{12}f_{D_3}(d)&=&\frac{d}{6}\left[-\frac{d^2}{3\sqrt{3}}\sin^{-1}\frac{
3\sqrt{3}}{2d}+\left(\frac{\pi}{9\sqrt{3}}-\frac{1}{3}\right)d^2+2d-\frac{\sqrt{
4d^2-27}}{4}-\frac{9}{4}\right], \nonumber
\end{eqnarray}
\begin{eqnarray}
\frac{1}{12}f_{D_6}(d)&=&\frac{d}{6}\left[\left(\frac{d^2}{3\sqrt{3}}-\frac{4}{
\sqrt{3}}\right)\sin^{-1}\frac{\sqrt{3}}{d}
-\left(\frac{d^2}{3\sqrt{3}}+\frac{1}{2\sqrt{3}}\right)\sin^{-1}\frac{\sqrt{3}}{
2d} \right. \nonumber\\
&&\left.-\left(\frac{2d^2}{3\sqrt{3}}+\frac{7\sqrt{3}}{2}\right)\sin^{-1}\frac{
3\sqrt{3}}{2d}+\left(\frac{1}{3}+\frac{2\pi}{9\sqrt{3}}\right)d^2-\frac{\sqrt{
4d^2-3}}{4}-\sqrt{d^2-3}\right. \nonumber\\
&&\left.-\frac{5}{3}\sqrt{4d^2-27}+\frac{11\pi}{2\sqrt{3}}+\frac{13}{4}
\right], \nonumber
\end{eqnarray}
\begin{eqnarray}
\frac{1}{24}f_{D_7}(d)&=&\frac{d}{12}\left[\frac{1}{\sqrt{3}}\left(\frac{2d^2}{3
}+1\right)\sin^{-1}\frac{\sqrt{3}}{2d}-\frac{4}{\sqrt{3}}\left(\frac{d^2}{3}
+2\right)\sin^{-1}\frac{\sqrt{3}}{d}-3\sqrt{3}\sin^{-1}\frac{3\sqrt{3}}{2d}
\right. \nonumber\\
&&\left.+\left(\frac{1}{3}+\frac{2\pi}{9\sqrt{3}}\right)d^2+\frac{\sqrt{4d^2-3}}
{2}-4\sqrt{d^2-3}-\sqrt{ 4d^2-27}+\frac{17\pi}{3\sqrt{3}}+\frac{9}{2}\right],
\nonumber
\end{eqnarray}
\begin{eqnarray}
\frac{1}{24}f_{D_8}(d)&=&\frac{d}{12}\left[\frac{2}{\sqrt{3}}\left(\frac{d^2}{3}
+8\right)\sin^{-1}\frac{\sqrt{3}}{d}+\frac{4}{\sqrt{3}}\left(\frac{d^2}{3}
+6\right)\sin^{-1}\frac{3\sqrt{3}}{2d}-\left(1+\frac{2\pi}{3\sqrt{3}}
\right)d^2\right. \nonumber\\
&&\left.+6\sqrt{d^2-3}+\frac{11}{3}\sqrt{4d^2-27}-\frac{40\pi}{3\sqrt{3}}
-11\right]. \nonumber
\end{eqnarray}
Thus,
\begin{eqnarray}
f_{2D}(d)&=&\frac{3}{8}f_{D_{\rm
H_A}}(d)+\frac{1}{12}\left[f_{D_4}(d)+f_{D_3}(d)+f_{D_6}(d)\right]+\frac{1}{24}
\left[f_{D_7}(d)+f_{D_8}(d)\right] \nonumber\\
&=&\frac{d}{6}\left[-\frac{4}{\sqrt{3}}\left(\frac{2(d/2)^2}{3}+1\right)\sin^{-1
}\frac{\sqrt{3}}{2(d/2)}+\frac{2\pi}{3\sqrt{3}}\left(\frac{d}{2}\right)^2-2\sqrt
{4\left(\frac{d}{2}\right)^2-3}+\frac{10\pi}{3\sqrt{3}}\right]\nonumber \\
&=&\frac{1}{2}f_{D_{\rm H_I}}(\frac{d}{2}). \nonumber
\end{eqnarray}

\item{$3\leq d \leq 2\sqrt{3}:$}
\begin{eqnarray}
\frac{3}{8}f_{D_{\rm
H_A}}(d)&=&\frac{d}{12}\left[\frac{2}{\sqrt{3}}\left(\frac{d^2}{3}
+12\right)\sin^{-1}\frac{3\sqrt{3}}{2d}-\frac{4}{\sqrt{3}}\left(\frac{d^2}{3}
+2\right)\sin^{-1}\frac{\sqrt{3}}{d} \right. \nonumber\\
&&\left.+\left(\frac{\pi}{9\sqrt{3}}-\frac{1}{3}\right)d^2-4\sqrt{d^2-3}+\frac{
19}{6}\sqrt{4d^2-27}-\frac{8\pi}{3\sqrt{3}}-\frac{9}{2}\right], \nonumber
\end{eqnarray}
\begin{eqnarray}
\frac{1}{12}f_{D_6}(d)&=&\frac{d}{6}\left[\left(\frac{d^2}{3\sqrt{3}}-\frac{4}{
\sqrt{3}}\right)\sin^{-1}\frac{\sqrt{3}}{d}-\left(\frac{d^2}{3\sqrt{3}}+\frac{1}
{2\sqrt{3}}\right)\sin^{-1}\frac{\sqrt{3}}{2d}+\frac{5\sqrt{3}}{2}\sin^{-1}\frac
{3\sqrt{3}}{2d} \right. \nonumber\\
&&\left.-\frac{\sqrt{4d^2-3}}{4}-\sqrt{d^2-3}+\frac{5}{6}\sqrt{4d^2-27}-\frac{
\pi}{2\sqrt{3}}-\frac{5}{4}\right], \nonumber
\end{eqnarray}
\begin{eqnarray}
\frac{1}{24}f_{D_7}(d)&=&\frac{d}{12}\left[\frac{1}{\sqrt{3}}\left(\frac{2d^2}{3
}+1\right)\sin^{-1}\frac{\sqrt{3}}{2d}-\frac{4}{\sqrt{3}}\left(\frac{d^2}{3}
+2\right)\sin^{-1}\frac{\sqrt{3}}{d} \right. \nonumber\\
&&\left.+\left(\frac{2d^2}{3\sqrt{3}}+3\sqrt{3}
\right)\sin^{-1}\frac{3\sqrt{3}}{2d}+\frac{\sqrt{4d^2-3}}{2}-4\sqrt{d^2-3}+\frac
{3}{2}\sqrt{4d^2-27}-\frac{\pi}{3\sqrt{3}}\right], \nonumber
\end{eqnarray}
\begin{eqnarray}
\frac{1}{24}f_{D_8}(d)&=&\frac{d}{12}\left[\frac{2}{\sqrt{3}}\left(\frac{d^2}{3}
+8\right)\sin^{-1}\frac{\sqrt{3}}{d}-\frac{4}{\sqrt{3}}\left(\frac{d^2}{3}
+12\right)\sin^{-1}\frac{3\sqrt{3}}{2d}+\left(\frac{1}{3}+\frac{2\pi}{9\sqrt{3}}
\right)d^2 \right. \nonumber\\
&&\left.+6\sqrt{d^2-3}-\frac{19}{3}\sqrt{4d^2-27}+\frac{32\pi}{3\sqrt{3}}
+7\right]. \nonumber
\end{eqnarray}
Thus,
\begin{eqnarray}
f_{2D}(d)&=&\frac{3}{8}f_{D_{\rm
H_A}}(d)+\frac{1}{12}f_{D_6}(d)+\frac{1}{24}\left[f_{D_7}(d)+f_{D_8}(d)\right]
\nonumber\\
&=&\frac{d}{6}\left[-\frac{4}{\sqrt{3}}\left(\frac{2(d/2)^2}{3}+1\right)\sin^{-1
}\frac{\sqrt{3}}{2(d/2)}+\frac{2\pi}{3\sqrt{3}}\left(\frac{d}{2}\right)^2-2\sqrt
{4\left(\frac{d}{2}\right)^2-3}+\frac{10\pi}{3\sqrt{3}}\right]\nonumber \\
&=&\frac{1}{2}f_{D_{\rm H_I}}(\frac{d}{2}). \nonumber
\end{eqnarray}

\item{$2\sqrt{3} \leq d\leq \sqrt{13}:$}
\begin{eqnarray}
\frac{3}{8}f_{D_{\rm
H_A}}(d)&=&\frac{d}{12}\left[\frac{2}{\sqrt{3}}\left(\frac{d^2}{3}
+12\right)\left(\sin^{-1}\frac{3\sqrt{3}}{2d}+\sin^{-1}\frac{2\sqrt{3}}{d}
\right)-\left(\frac{2}{3}+\frac{4\pi}{9\sqrt{3}}\right)d^2\right. \nonumber\\
&&\left.+\frac{19}{6}\sqrt{4d^2-27}+\frac{16}{3}\sqrt{d^2-12}-\frac{16\pi}{\sqrt
{3}}-\frac{25}{2}\right], \nonumber
\end{eqnarray}
\begin{eqnarray}
\frac{1}{12}f_{D_6}(d)&=&\frac{d}{6}\left[\left(\frac{d^2}{3\sqrt{3}}+\frac{8}{
\sqrt{3}}\right)\sin^{-1}\frac{2\sqrt{3}}{d}-\left(\frac{d^2}{3\sqrt{3}}+\frac{1
}{2\sqrt{3}}\right)\sin^{-1}\frac{\sqrt{3}}{2d}+\frac{5\sqrt{3}}{2}\sin^{-1}
\frac{3\sqrt{3}}{2d}\right. \nonumber\\
&&\left.-\left(\frac{1}{6}+\frac{\pi}{9\sqrt{3}}\right)d^2-\frac{\sqrt{4d^2-3}}{
4}+\frac{5}{6}\sqrt{4d^2-27}+2\sqrt{d^2-12}-\frac{31\pi}{6\sqrt{3}}-\frac{9}{4}
\right], \nonumber
\end{eqnarray}
\begin{eqnarray}
\frac{1}{24}f_{D_7}(d)&=&\frac{d}{12}\left[\frac{1}{\sqrt{3}}\left(\frac{2d^2}{3
}+1\right)\sin^{-1}\frac{\sqrt{3}}{2d}+\left(\frac{2d^2}{3\sqrt{3}}+3\sqrt{3}
\right)\sin^{-1}\frac{3\sqrt{3}}{2d}+\frac{8}{\sqrt{3}}\sin^{-1}\frac{2\sqrt{3}}
{d}\right. \nonumber\\
&&\left.-\left(\frac{1}{3}+\frac{2\pi}{9\sqrt{3}}\right)d^2+\frac{\sqrt{4d^2-3}}
{2}+\frac{3}{2}\sqrt{4d^2-27}+\frac{4}{3}\sqrt{d^2-12}-\frac{17\pi}{3\sqrt{3}}
-8\right], \nonumber
\end{eqnarray}
\begin{eqnarray}
\frac{1}{24}f_{D_8}(d)&=&\frac{d}{12}\left[-\frac{4}{\sqrt{3}}\left(\frac{d^2}{3
}+12\right)\sin^{-1}\frac{3\sqrt{3}}{2d}-\frac{2}{\sqrt{3}}\left(\frac{d^2}{3}
+8\right)\sin^{-1}\frac{2\sqrt{3}}{d}+\left(1+\frac{2\pi}{3\sqrt{3}}\right)d^2
\right. \nonumber\\
&&\left.-\frac{19}{3}\sqrt{4d^2-27}-4\sqrt{d^2-12}+\frac{64\pi}{3\sqrt{3}}
+17\right]. \nonumber
\end{eqnarray}
Thus,
\begin{eqnarray}
f_{2D}(d)&=&\frac{3}{8}f_{D_{\rm
H_A}}(d)+\frac{1}{12}f_{D_6}(d)+\frac{1}{24}\left[f_{D_7}(d)+f_{D_8}(d)\right]
\nonumber\\
&=&\frac{d}{6}\left[\frac{1}{\sqrt{3}}\left(\frac{d^2}{3}+16\right)\sin^{-1}
\frac{2\sqrt{3}}{d}-\left(\frac{\pi}{9\sqrt{3}}+\frac{1}{6}\right)d^2+\frac{10}{
3}\sqrt{d^2-12}-\frac{16\pi}{3\sqrt{3}}-4\right]\nonumber\\
&=&\frac{d}{6}\left[\frac{4}{\sqrt{3}}\left(\frac{(d/2)^2}{3}+4\right)\sin^{-1}
\frac{\sqrt{3}}{d/2}-\left(\frac{4\pi}{9\sqrt{3}}+\frac{2}{3}\right)\left(\frac{
d}{2}\right)^2+\frac {20}{3}\sqrt{\left(\frac{d}{2}\right)^2-3}
-\frac{16\pi}{3\sqrt{3}}-4 \right]\nonumber\\
&=&\frac{1}{2}f_{D_{\rm H_I}}(\frac{d}{2}). \nonumber
\end{eqnarray}

\item{$\sqrt{13} \leq d\leq 4:$}
\begin{eqnarray}
\frac{1}{24}f_{D_8}(d)=\frac{d}{12}\left[\frac{2}{\sqrt{3}}\left(\frac{d^2}{3}
+16\right)\sin^{-1}\frac{2\sqrt{3}}{d}-\left(\frac{1}{3}+\frac{2\pi}{9\sqrt{3}}
\right)d^2+\frac{20}{3}\sqrt{d^2-12}-\frac{32\pi}{3\sqrt{3}}-8\right]. \nonumber
\end{eqnarray}
Thus,
\begin{eqnarray}
f_{2D}(d)&=&\frac{1}{24}f_{D_8}(d)=\frac{d}{6}\left[\frac{1}{\sqrt{3}}
\left(\frac{d^2}
{3}+16\right)\sin^{-1}\frac{2\sqrt{3}}{d}-\left(\frac{\pi}{9\sqrt{3}}+\frac{1}{6
}\right)d^2+\frac{10}{3}\sqrt{d^2-12}-\frac{16\pi}{3\sqrt{3}}-4\right]
\nonumber\\
&=&\frac{d}{6}\left[\frac{4}{\sqrt{3}}\left(\frac{(d/2)^2}{3}+4\right)\sin^{-1}
\frac{\sqrt{3}}{d/2}-\left(\frac{4\pi}{9\sqrt{3}}+\frac{2}{3}\right)\left(\frac{
d}{2}\right)^2+\frac{20}{3}\sqrt{\left(\frac{d}{2}\right)^2-3}-\frac{16\pi}{
3\sqrt{3}}-4\right]\nonumber\\
&=&\frac{1}{2}f_{D_{\rm H_I}}(\frac{d}{2}). \nonumber
\end{eqnarray}
\end{enumerate}

In summary, we have $f_{2D}(d)=\frac{1}{2}f_{D_{\rm H_I}}(\frac{d}{2})$ by recursion, and
the probabilistic sum $\frac{3}{16}f_{D_{\rm H_I}}(d)+\frac{3}{8}f_{D_{\rm
H_A}}(d)+\frac{1}{12}\left[f_{D_2}(d)+f_{D_4}(d)+f_{D_3}(d)+f_{D_6}(d)\right]
+\frac{1}{24}\left[f_{D_7}(d)+f_{D_8}(d)\right]+\frac{1}{48}f_{D_{\rm
R_I}}(d)$ is equal to $\frac{1}{2}f_{D_{\rm H_I}}(\frac{d}{2})$ in all the cases discussed
above. The results are a strong validation of the correctness of the distance 
distribution functions that we have derived.

\section{Practical Results}

\subsection{Statistical Moments of Random Distances}

The distance distribution functions given in Section~\ref{sec:result} can conveniently lead to all the
statistical moments of the random distances associated with hexagons. Given $f_{D_{\rm H_I}}(d)$ in
(\ref{eq:fd_h_within}), for example, the first moment (mean) of $d$, or the average distance within a regular 
hexagon, is
\begin{eqnarray}
 M_{D_{\rm H_I}}^{(1)}=\int_0^2xf_{D_{\rm
H_I}}(x)dx=\frac{7\sqrt{3}}{30}-\frac{7}{90}+\frac{1}{60}\left[28\ln
\left(2\sqrt{3}+3\right)+29\ln \left(2\sqrt{3}-3\right)\right]\approx
0.82625\textcolor{red}{89495}, \nonumber
\end{eqnarray}
and the second raw moment is
\begin{eqnarray}
 M_{D_{\rm H_I}}^{(2)}=\int_0^2x^2f_{D_{\rm H_I}}(x)dx=\frac{5}{6},
\nonumber
\end{eqnarray}
from which the variance, or the second central moment, can be derived as
\begin{eqnarray}
 Var_{D_{\rm H_I}}=M_{D_{\rm H_I}}^{(2)}-\left[M_{D_{\rm H_I}}^{(1)}\right]^2\approx
0.150629\textcolor{red}{4817}. \nonumber
\end{eqnarray}

When the side length of a hexagon is scaled by $s$, the corresponding first 
two statistical moments given above then become
\begin{equation}
\label{within-moments-numerical}
M_{D_{\rm H_I}}^{(1)}=0.82625\textcolor{red}{89495}s,~~\mbox{}~~M_{D_{\rm H_I}}^{(2)}=\frac{5s}{6}
~~\mbox{ and }~~Var_{D_{\rm H_I}}=0.150629\textcolor{red}{4817}s^2.
\end{equation}

\begin{table}
  \caption{Moments and Variance---Numerical vs Simulation Results}
  \centering
  \begin{tabular}{|c||c|c|c|c|}
    \hline
    Geometry & PDF/Sim & $M_{D}^{(1)}$ & $M_{D}^{(2)}$ & $Var_{D}$ \\ \hline 
\hline 
   Within a &  $f_{D_{\rm H_I}}(d) $ & $0.82625\textcolor{red}{89495}s$ & $0.8333333333s$ &
$0.150629\textcolor{red}{4817}s^2$\\ 
    \cline{2-5}
    Single Hexagon & Sim & $ 0.8263306317s$ & $0.8335924725s$ & $0.1507701596s^2$ 
\\ \hline
    Between two & $f_{D_{\rm H_A}}(d)$ & $1.8564318344s$ & $3.832947195s$ &
$0.3866080394s^2$\\ 
    \cline{2-5}
    Adjacent Hexagons & Sim & $1.8583366966s$ & $3.8326819696s$ & $0.3792666917s^2$
\\ \hline
  \end{tabular}
  \label{tab:moment}
\end{table}

Table~\ref{tab:moment} lists the first two moments, and the variance of the random
distances in the two cases given in Section~\ref{sec:result}, and the
corresponding simulation results for verification purposes.

\subsection{Polynomial Fits of Random Distances}

\begin{table}
  \caption{Coefficients of the Polynomial Fit and the Norm of Residuals (NR)}
  \centering
  \begin{tabular}{|c||c|c|c|}
    \hline
    PDF & Degree & Polynomial Coefficients & NR \\ \hline \hline
    & & $10^2\times
\left[-0.0146710~~0.136604~-0.538052~~1.167903~\right.$ & \\
    $f_{D_{\rm H_I}}(d)$ & $10$ &
$-1.525478~~1.230615~-0.605940~~0.175147~-0.043772$ & $0.075608$ \\
    & & 
$\left.0.025830~-0.000025\right]$ & \\ \hline
    & & $10^4\times
\left[0.00000035~-0.000013~~0.000207~-0.002094\right.$ & \\
    & &
$0.014469~-0.072522~~0.272508~-0.782682~~1.736254$ & \\
    $f_{D_{\rm H_A}}(d)$ & $20$ & $-2.986406~~3.976655~-4.072372
~~3.169347~-1.841066$ & $0.191157$ \\
    & & $0.778001~-0.230634~~0.045522-0.005534~~0.000394 $ & \\
    & & $\left.-0.0000103~~0.00000007092\right]$ & \\\hline
  \end{tabular}
  \label{tab:poly}
\end{table}

\begin{figure}
\centering
  \subfloat[Within a Single
Hexagon]{\includegraphics[width=0.5\columnwidth]{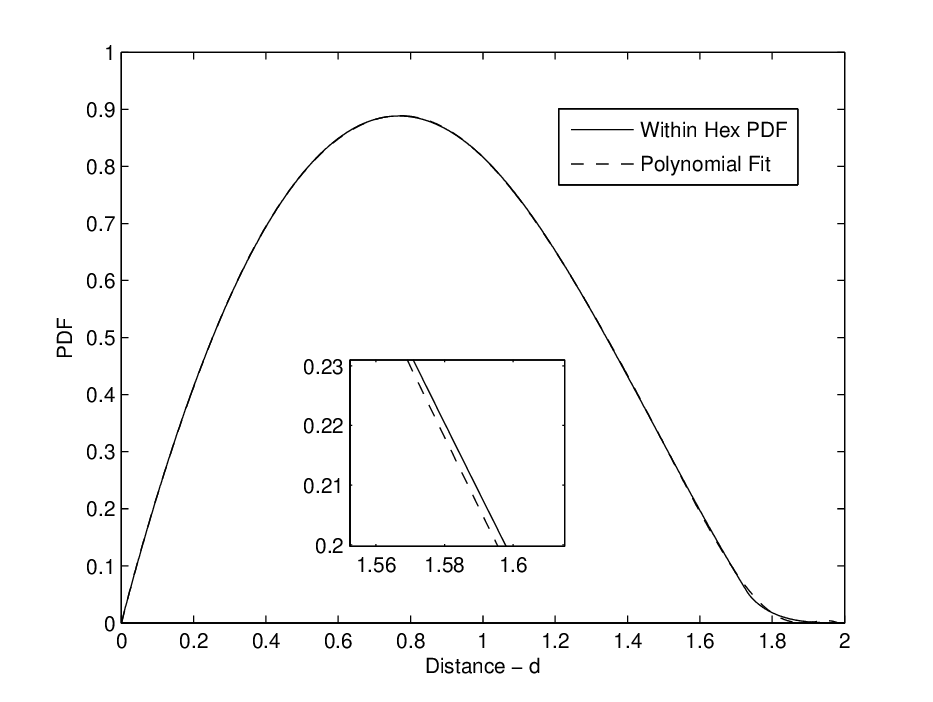}}
  \subfloat[Between two Adjacent
Hexagons]{\includegraphics[width=0.5\columnwidth]{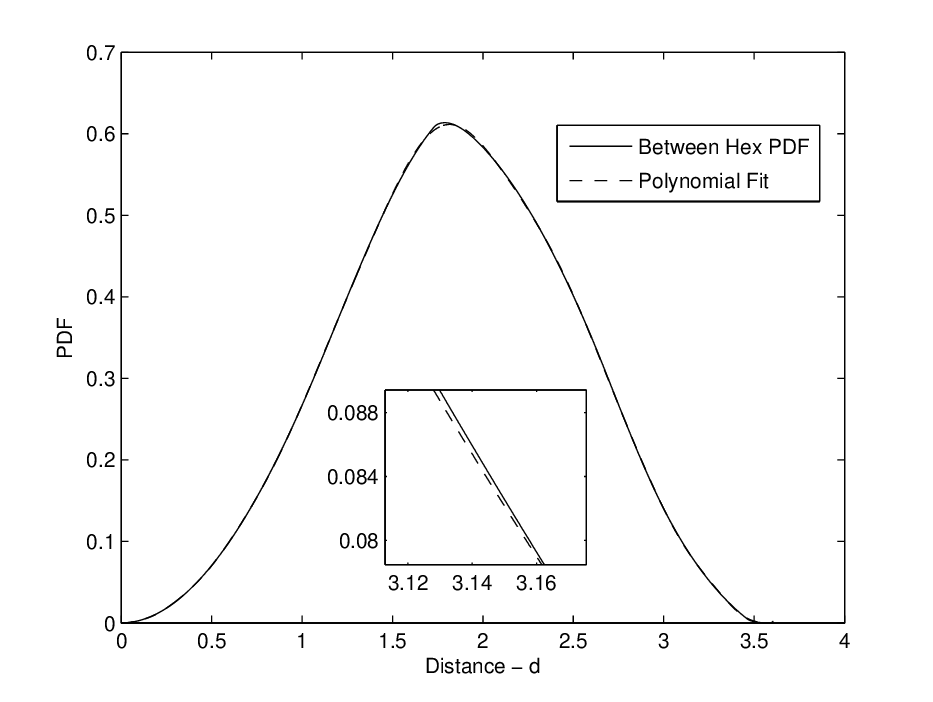}}
  \caption{Polynomial Fit.}
  \label{fig:hex_poly}
\end{figure}

Table~\ref{tab:poly} lists the coefficients of the high-order polynomial fits of the original 
PDFs given in Section~\ref{sec:result}, from the highest degree (degree-$10$ 
for (\ref{eq:fd_h_within}) and degree-$20$ for (\ref{eq:fd_h_btw})) to $d^0$, together with the corresponding norm of residuals.
Figure~\ref{fig:hex_poly} (a)--(b) plot the polynomials listed in
Table~\ref{tab:poly} with the original PDFs. From the figure, it can be seen that both polynomials match closely
with the original PDFs. These high-order polynomials facilitate further manipulations of the distance distribution functions, with a high accuracy.

\section{Conclusions}
\label{sec:conclude}

In this report, we gave the closed-form probability density
functions of the random distances associated with hexagons. The correctness of the
obtained results has been verified by a recursion and a probabilistic sum,
in addition to simulation. The first two statistical moments, and the 
polynomial fits of the density functions are also given for practical uses.

\section*{Acknowledgment}

This work is supported in part by the NSERC, CFI and BCKDF. The authors also
want to thank Dr. Lin Cai for the initial work involving squares, Dr. Aaron
Gulliver for posing the new problem associated with hexagons, and Dr. A.M.
Mathai, Dr. David Chu and Dr. T.S. Hale for their correspondence and
encouragement.

\textcolor{red}{In 2021, we would like to thank Dr.~Uwe B\"asel~\cite{basel21moments} for pointing out a numerical calculation error of the values of $M_{D_{\rm H_I}}^{(1)}$ in (\ref{within-moments-numerical}) and thus $Var_{D_{\rm H_I}}$ (i.e., the first moment and variance of the random distances within a regular hexagon of side length $s$) in the 2011 version of this report, which now have been corrected in red on page~\pageref{within-moments-numerical}. The original expression of $M_{D_{\rm H_I}}^{(1)}$, however, as confirmed by Dr.~B\"asel, is correct and equivalent to $M_1$ in \cite{basel21moments}, as
\begin{equation}
\begin{array}{rl}
\frac{M_{D_{\rm H_I}}^{(1)}}{s}=&\frac{7\sqrt{3}}{30}-\frac{7}{90}+\frac{1}{60}[28\ln(2\sqrt{3}+3)+29\ln(2\sqrt{3}-3)]\\
=&-\frac{7}{90}+\frac{7}{10\sqrt{3}}+\frac{1}{60}[29\ln(2\sqrt{3}+3)-\ln(2\sqrt{3}+3)+29\ln(2\sqrt{3}-3)]\\
=&-\frac{7}{90}+\frac{7}{10\sqrt{3}}+\frac{1}{60}\{29[\ln(2\sqrt{3}+3)+\ln(2\sqrt{3}-3)]-[\ln(2+\sqrt{3})+\ln{\sqrt{3}}]\}\\
=&-\frac{7}{90}+\frac{7}{10\sqrt{3}}+\frac{1}{60}\{29\ln[(2\sqrt{3}+3)(2\sqrt{3}-3)]-[\ln(2+\sqrt{3})+\ln{\sqrt{3}}]\}\\
=&-\frac{7}{90}+\frac{7}{10\sqrt{3}}+\frac{1}{60}[29\ln{3}-\ln(2+\sqrt{3})-\frac{1}{2}\ln{3}]\\
=&-\frac{7}{90}+\frac{7}{10\sqrt{3}}+(\frac{29}{60}-\frac{1}{120})\ln{3}-\frac{\ln(2+\sqrt{3})}{60}\\
=&-\frac{7}{90}+\frac{7}{10\sqrt{3}}+\frac{57}{120}\ln{3}-\frac{\ln(2+\sqrt{3})}{60}\\
=&-\frac{7}{90}+\frac{7}{10\sqrt{3}}+\frac{19\ln{3}}{40}-\frac{\ln(2+\sqrt{3})}{60}\\
=&\frac{M_1}{\ell_1},
\end{array}
\end{equation}
where $M_1$ is defined in \cite{basel21moments} as the first moment of the random distances in a regular hexagon of side length $\ell_1$, $r$ is the radius of the hexagon's circumscribed circle, and $\ell_1=r=s$ by definition. This, in turn, provides an independent verification of the results in \cite{rhombus} and this report.}

\end{document}